# Giga-Periodic Orbits for Weakly Coupled Tent and Logistic Discretized Maps

## René Lozi

Laboratoire J.A. Dieudonné – UMR du CNRS N° 6621
University of Nice-Sophia-Antipolis, Parc Valrose, 06108 Nice Cedex 02, France

and

Institut Universitaire de Formation des Maîtres de l'académie de Nice,
89 avenue George V, 06046 Nice Cedex 1, France
rlozi@unice.fr

___

Abstract

Simple dynamical systems often involve periodic motion. Quasi-periodic or chaotic motion is frequently present in more complicated dynamical systems. However, for the most part, underneath periodic motion models chaotic motion. Chaotic attractors are nearly always present in such dissipative systems. Since their discovery in 1963 by E. Lorenz, they have been extensively studied in order to understand their nature. In the past decade, the aim of the research has been shifted to the applications for industrial mathematics. Their importance in this field is rapidly growing.

Chaotic orbits embedded in chaotic attractor can be controlled allowing the possibility to control laser beams or chemical processes and improving techniques of communications. They can also produce very long sequences of numbers which can be used as efficiently as random numbers even if they have not the same nature.

However, mathematical results concerning chaotic orbits are often obtained using sets of real numbers (belonging to **R** or **R**$^n$) (*e.g.* the famous theorem of A. N. Sharkovskiĭ which defines which ones periods exist for continuous functions such as logistic or tent maps).

O.E. Lanford III reports the results of some computer experiments on the orbit structure of the discrete maps on a finite set which arise when an expanding map of the circle is iterated "naively" on the computer. There is a huge gap between these results and the theorem of Sharkovskiĭ, due to the discrete nature of floating points used by computers.

This article introduces *new models of very very weakly coupled logistic and tent maps* for which orbits of very long period are found. The length of these periods is far greater than one billion. We call giga-periodic orbits such orbits for which the length is greater than $10^9$ and less than $10^{12}$. Tera, and peta periodic orbits are the name of the orbits the length of which is one thousand or one million greater. The property of these models relatively to the distribution of the iterates (invariant measure) are described. They are found very useful for industrial mathematics for a variety of purposes such as generation of cryptographic keys, computer games and some classes of scientific experiments.

___

**Key words :** chaotic numbers, coupled tent maps, coupled logistic maps, giga-periodic orbit, Lozi map, Henon map.

## 1. Introduction

Randomness and random numbers have traditionally been used for a variety of purposes, for example games such as dice games on computer or generation of cryptographic keys. With the advent of computers, people recognized the need for a means of introducing ramdomness into a computer program. Surprisingly, it is difficult to get a computer to do something by chance. A computer running a program follows its instructions blindly and is therefore completely predictable.

Computer engineers chose to introduce randomness into computers in the form of *pseudo-random number generators* (PRNG). Pseudo-random numbers are not truly random and are rather computed from a mathematical formula or simply taken from a precalculed list. A lot of research has gone into pseudo-random number theory and modern algorithms for generating them are so good that the numbers look like they were really random. In fact, we never prove that a string of number is random (such as the decimals of $\pi$), but only that it is not random by finding some pattern in the numbers, and with so many possible patterns, testing them all is a hopeless undertaking. Pseudo-random numbers have the characteristic that they are predictable, i.e. they can be predicted if we know in the sequence the first number taken from (called the *seed* or the *root*). For some purposes, predictabilty is a good characteristic, for others it is not.

Physicists believe that atomic processes governed by the laws of quantum machanics are inherently probabilistic and random. *True random number generators* based on such processes have been built but often perform poorly compared with deterministic PRNG [1].

Chaos is used to describe fundamental disorder generated by simple deterministic systems with only a few elements [2]. However, chaos is not strictly equivalent to random even if the irregularities of chaotic and random sequences in the time domain are often quite similar. As an illustration, the logistic sequence and random sequence are plotted as shown in Fig. 1 (a, b) where it is difficult to observe any difference. However, by plotting the phase space of chaotic and random sequences, (Fig. 2 (a, b)), the chaotic sequence can be easily identified because there is a regular pattern (a parabola) in the phase space plot. Consecutive points of a chaotic sequence are highly correlated, but for the case of a pure random sequence.

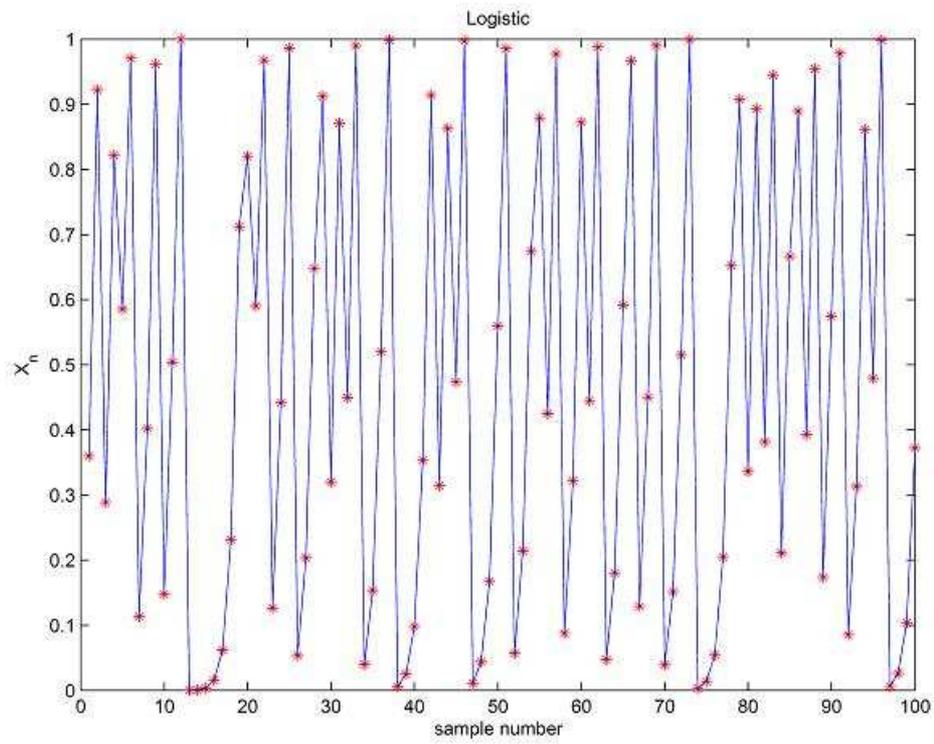

(a)

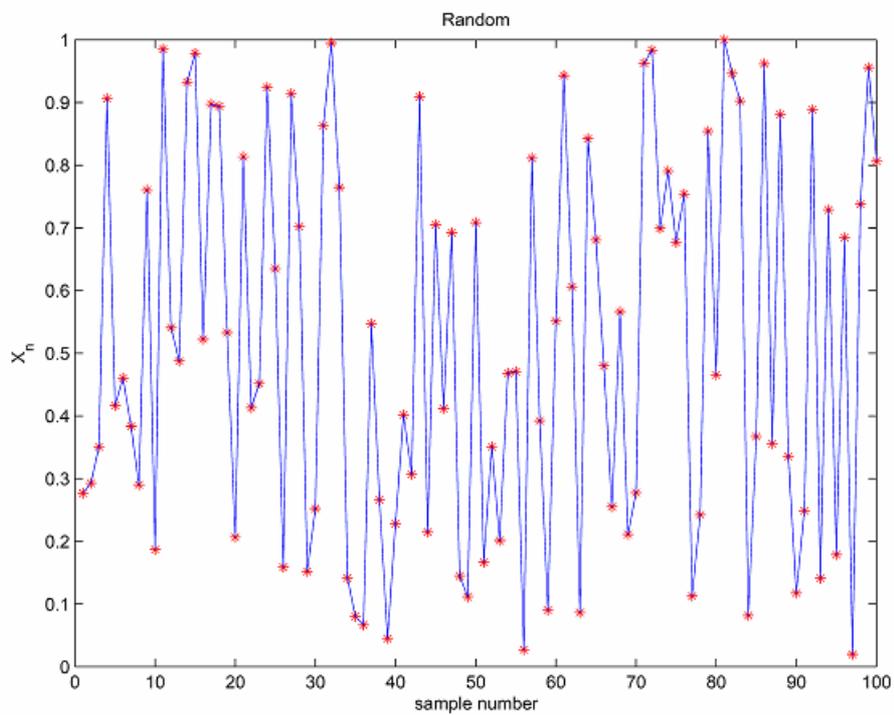

(b)

**Fig. 1. (a) Logistic and (b) random sequences, respectively, with** $a = 4$ **(from [4]).**

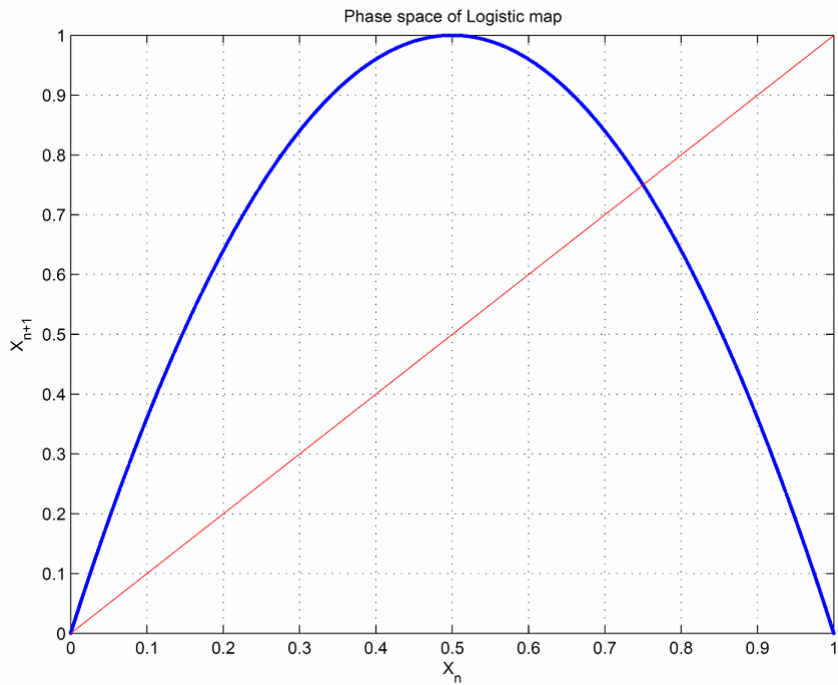

(a)

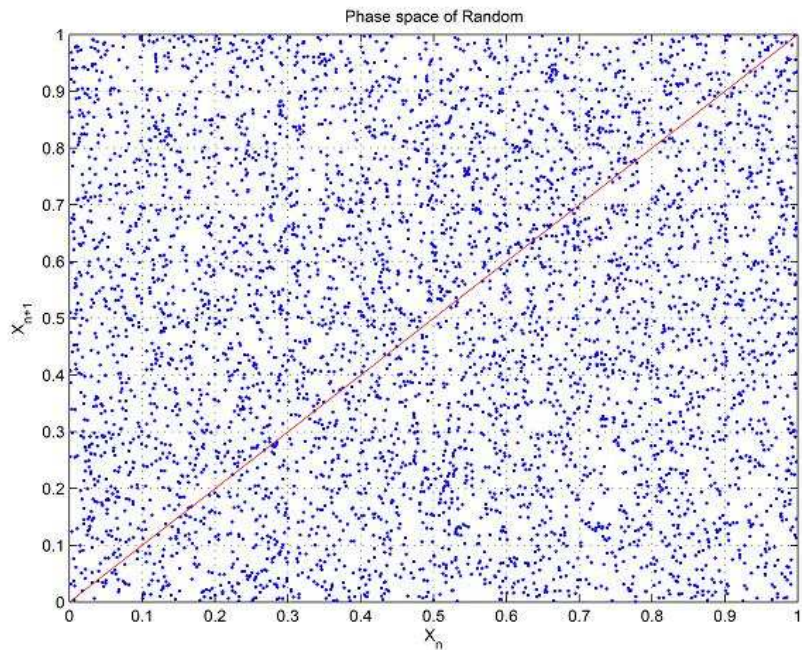

(b)

**Fig. 2. Phase space plot of: (a) the logistic sequence with *a* = 4 and (b) random sequence (from [4]).**

The difference between chaos and random has been recently carrefuly studied [3]. For several applications chaos has been found more useful that random in some cases [4].

It is noteworthy that the new models of very weakly coupled maps introduced in this article are more powerful than the usual chaotic maps used to generate chaotic sequences.

## 2. Dynamical Systems and Chaos

Dynamical systems which model many physical laws have been extensively studied for more than three centuries (at least since the Newton's idea of modeling the motion of physical systems with equations). The possible existence of chaotic dynamicals systems has been pointed out by Poincaré [5] who infered that the three-body problem (i.e. the gravitational interaction between sun, moon and earth) was impossible to solve.

Now we explain what is a dynamical systems.

The function

$$f(x) = 2.4\, x \qquad (2.1)$$

is a rule that assigns to each number $x$ a number 2.4 times as large. This is a simple mathematical model. We may imagine that $x$ be the range of the population of tame rabitts in a rabbit-hutch and $f(x)$ denotes the range of the same population two months later. Then the rule expresses the fact that the population is multiplied by 2.4 every two months and can be represented as

$$x_{n+1} = 2.4\, x_n \qquad (2.2)$$

where $x_n$ is for the population of tame rabitts for a given month and $x_{n+1}$ stands for the population two months later.

A *dynamical system* consists of a set of possible states, together with a rule that determines the present state in terms of past states. In the previous paragraph, we discussed a simple dynamical system whose states are population levels, that change with time under the rule $x_{n+1} = f(x_n) = 2.4\, x_n$, where $n$ is the time and $x_n$ the population at time $n$. We will require that the rule is *deterministic*, which means that we can determine the present state *uniquely* from the past states.

Fundamental to science is the presumption that experiments are predictable and repeatable. Thus, it surprised most scientists when simple deterministic systems were found to be neither predictable nor repeatable. Instead they exhibited *chaos*, in which the tiniest change in the initial conditions produces a very different outcome, even when the governing equations are known exactly. Chaos occurs in several fields of physics: the solar system which seems perfectly regular, exhibiting very regular astronomical motions of the planets orbiting around the sun, is in fact chaotic on some aspects. Hyperion[1] (Fig. 3), one of the smaller moons of Saturn, is irregularly shaped and tumbles chaotically as its orbits Saturn [6].

---

[1] Hyperion, the sixteenth moon of Saturn, was discovered in 1848 by William Cranch Bond, his son George Phillips and independently by William Lassell one day after.

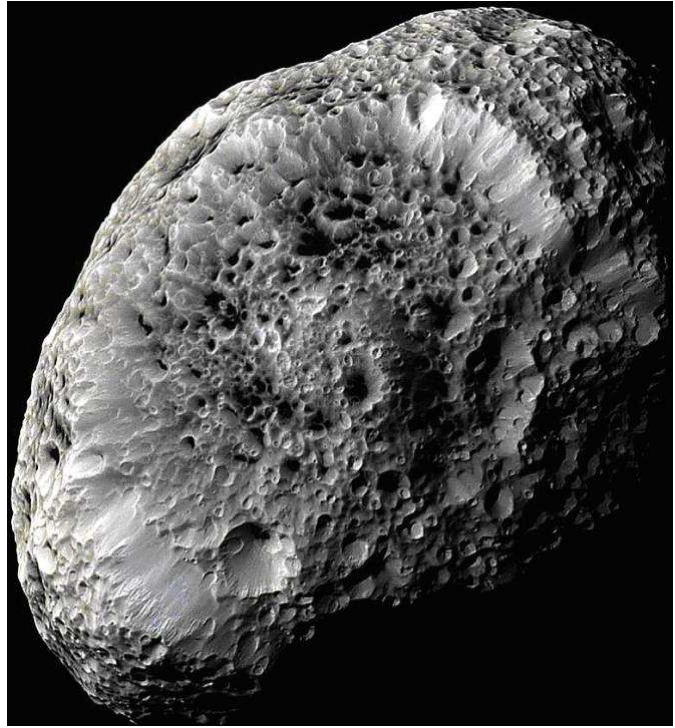

**Fig. 3. Image of a very spongy-looking Hyperion captured by *Cassini* (NASA/JPL). Hyperion is an irregularly shaped satellite, approximately 370 x 280 x 226 km in dimension. Hyperion orbits Saturn at a distance of almost 1.5 million km and takes over 3 weeks to complete one orbit.**

Halley's Comet whose orbit is highly elliptical is expected to be chaotic because its different interactions with the massive planets, such as Jupiter, on successive orbits [7]. Concerning the planets, their respective interactions also induce some effect on the tilt (obliquity) of their equatorial plane over their orbital plane. Mars presents a very large chaotic zone for its obliquity ranging (Fig. 4) from 0 to more than 60 degrees [8]. Two possible solutions for the past evolution of the obliquity of Mars, obtained using advanced numerical integration of the orbital motion of all the planets, and initial conditions and parameter within the uncertainty of the best known values for Mars rotational parameters are plotted on fig. 4 [9].

Fluids dynamics is also very rich in chaotic motions. A fluid (liquid, gas or plasma) consisting of very many infinitesimal fluid elements moves in response to their neighbours. The mixing that occurs while stirring the cream in a coffee cup is a chaotic process. Two sticks placed in a smoothly flowing river will remain side by side for a long time, but they will quickly separate if the river is turbulent.

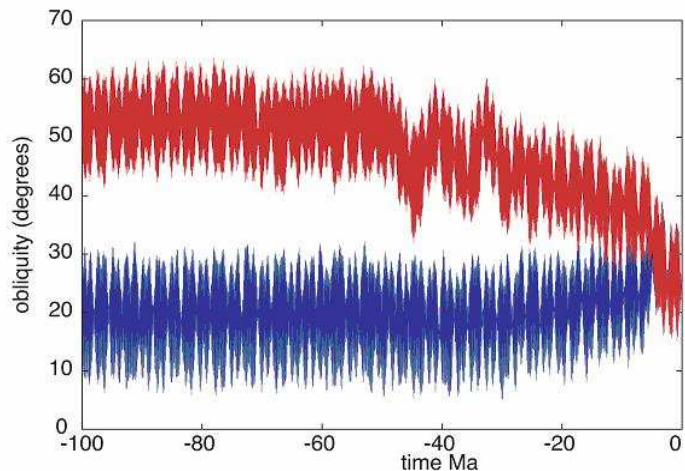

**Fig. 4.** Two different solutions for Mars obliquity over 100 Myrs, obtained with direct integration of the planetary orbits, and with initial conditions of the planet spin within the uncertainty of the most recent determinations(from [9]).

Chaotic atmospheric motion prevents long-term weather prediction and illustrates the sensitive dependence on initial conditions known as the *butterfly effect* in which a butterfly flapping its wings in south America can cause storms elsewhere in the world even in Europa [10]. Several other examples together with a very good description of chaotic phenomena and the mathematical treatment which can be done to the time-series issued from record of physical variable of chaotic system are described by Sprott [11].

We will emphasize two types of dynamical systems viz, *discrete-* and *continuous-time* systems. A discrete-time system takes the current time as input and updates the situation by producing a new state as output. A continuous-time is essentially the limit of discrete systems with smaller and smaller updating times. The governing rule in that case becomes a set of differential equations. A good introduction of mathematical tools actually used in order to study both systems can be found in Alligood et al. [12].

Scientists study models because they suggest how real-world processes behave. Every model of a physical process is at best an idealization. The goal of the model is to capture some feature of the physical process.

Generally, as mathematical tools are actually not very efficient in this beginning of 21st century, scientists are using chain of models instead of isolated ones. In such a chain of models, each one is consistent with the previous one. There is a simplification of the mathematical tools needed to study the behavior of the models from the beginning to the end of the chain:

Physical process in the real-world → Mathematical model using Partial Differential Equation (P.D.E.) → Mathematical model using Ordinary Differential Equation (O.D.E.) → Mathematical model involving Mapping (discretization) [13].

Hunt et al. [14] collect a coherent collection of articles related to the topics of "chaotic attractors" and the corresponding "natural invariant measures" that reflect the theoretical development of chaotic attractors. They are interesting to be read, and are appropriate for a graduate student seminar.

For in-depth study of chaotic dynamical systems refer to [12, 13, 14].

## 3. Mega and Giga-periodic Orbits for One- and Two-dimensional Chaotic Generators.

### 3.1. The Logistic Map, Cobweb Diagram and Periodic Orbit

*(a) The Logistic Map*

In previous section we considered the function (2.1) $f(x) = 2.4x$ or more generally

$$f_a(x) = a\,x \qquad (3.1)$$

as a very simple example of the rule which governs a dynamical system modelling the growing of tame rabitts population (*a* belonging to **R**). This linear function leads to exponential growth if $a > 1$. However exponential growth of rabitts population cannot continue for ever (even in Australia[2]). Hence linear dynamical system is not a realistic model of any natural process. Typically some nonlinearity stops or even reverses the growth[3]. The nonlinearity is negligible at small values of $x$ but dominates as $x$ increases. The simplest way to modify the function (3.1) to include a term that reduce the growth as $x$ increases is to substract one quadratic term $ax^2$.

The very interesting function obtained by this mean

$$f_a(x) = a\,x(1 - x) \qquad (3.2)$$

is the *logistic map*. The dynamical systems is defined by the *logistic equation* [11, 12]

$$x_{n+1} = f_a(x_n) = a\,x_n(1 - x_n) \qquad (3.3)$$

*(b) Cobweb Diagram*

A graph of function (3.2), called *logistic function* or *logistic curve*, is a parabola, as in Fig. 5. which shows that the case $a = 4$, which has the special property that it maps the unit interval ($0 \leq x \leq 1$) onto itself twice (once for $0 \leq x < 0.5$ and again for $0.5 < x \leq 1$ in the reverse direction). In the rubber-band analogy, this case corresponds to streching the rubber band to twice its lenght with each *iteration*, albeit nonuniformily and then folding it back onto itself. The map is one-dimensional since there is a single variable $x$ and the resulting curve is a line, even when it is plotted in a two-dimensional space.

---

[2] Wild rabbits have been causing widespread devastation in Australia for over 100 years. They were originally introduced to Australia by the First Fleet in 1788, but the current major infestation appears to be the results of 24 wild rabbits released by Thomas Austin on his Barwon Heads property in 1859 for hunting purposes. Many other farms released their rabbits into the wild after Austin. Within ten years of the 1859 introduction, the 24 rabbits had multiplied so much that 2 million a year could be shot or trapped without leaving any noticeable effect on the population size (following approximately the linear dynamical systems (3.1) with $a =1.3$).

[3] Releasing rabbit-borne diseases has proven somewhat successful in controlling the population of rabbits in Australia. In 1950, *Myxomatosis* was released into the rabbit population which caused the rabbit population to drop from an estimated 600 million to around 100 million (following then a nonlinear dynamical system). Genetic resistance in the remaining rabbits allowed the population to recover 200-300 million by 1991.

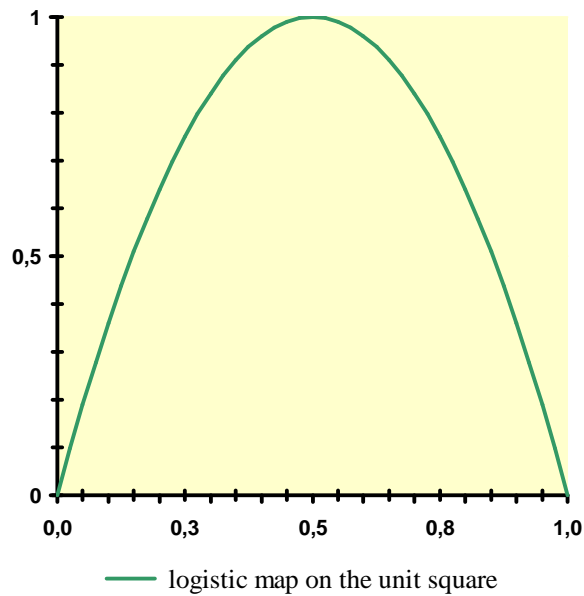

**Fig. 5.** **Graph of the logistic map on the unit square [0,1] x [0,1],** $f(x) = 4\,x\,(1-x)$.

In order to understand which states of the dynamical systems are obtained starting an initial value $x_0$ we ought to examine the sequence, $x_1, ..., x_n, x_{n+1}, ...$ generated by

$$\begin{aligned}
x_1 &= f_a(x_0) \\
x_2 &= f_a(x_1) = f_a \circ f_a(x_0) = f_a^{(2)}(x_0) \\
x_3 &= f_a(x_2) = f_a \circ f_a \circ f_a(x_0) = f_a^{(3)}(x_0) \\
&\vdots \\
x_{n+1} &= f_a(x_n) = \underbrace{f_a \circ \cdots \circ f_a \circ f_a}_{\ldots n+1 \text{ times} \ldots}(x_0) = f_a^{(n+1)}(x_0) \\
&\vdots
\end{aligned} \qquad (3.4)$$

Figure 6 shows a simple graphical way to do this. Start from the initial value $x_0$ on the horizontal axis. Draw a vertical line to the parabola to determine $x_1$. Then from that point draw a horizontal line to the line $x_{n+1} = x_n$. Then draw a vertical line to the parabola to determine $x_2$. Repeat as necessary to determine $x_n, x_{n+1}, \ldots$ Such a diagram is called a cobweb diagram for obvious reasons.

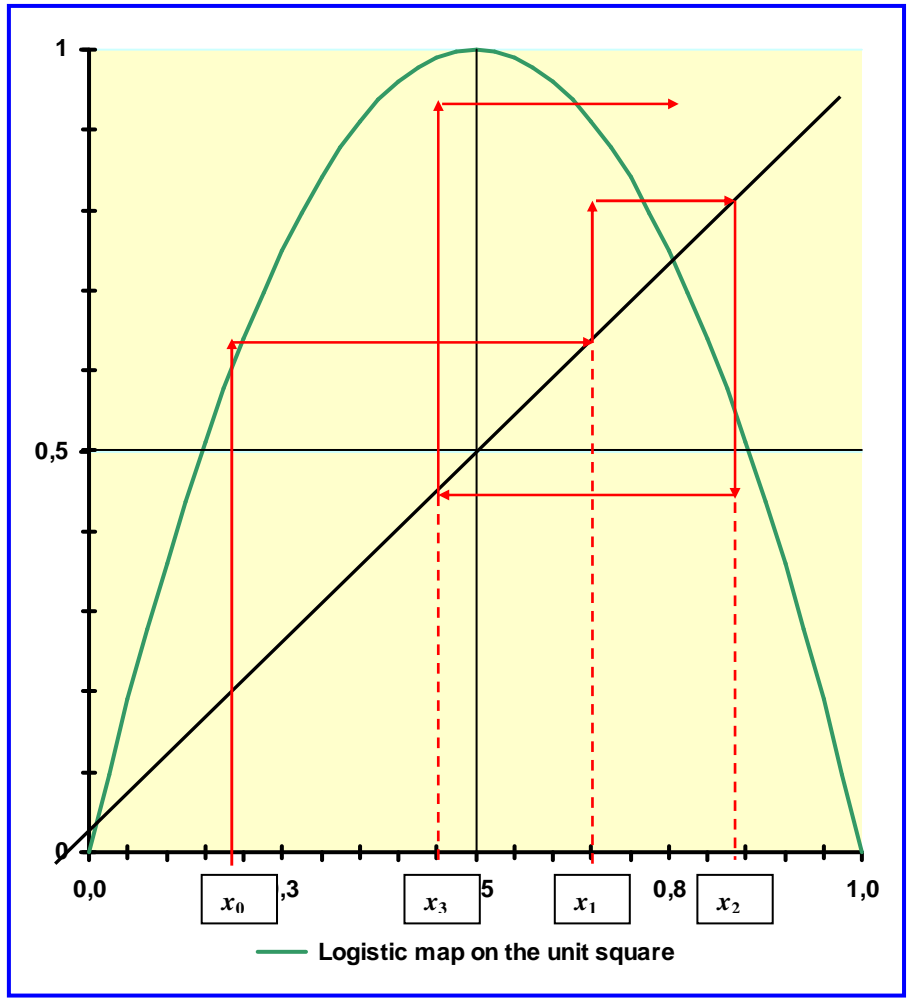

**Fig. 6.** Cobweb diagram for the logistic map on the unit square [0,1] x [0,1], $f(x) = 4\,x\,(1 - x)$.

*(c) Fixed Points and Periodic Orbits*

Consider Eq. (3.4) i.e. $\quad x_{n+1} = f_a(x_n) = f_a^{(n+1)}(x_0)\quad$ with $n = 0$.

If $x_0$ verifies

$$x_0 = f_a^1(x_0) = f_a(x_0) \tag{3.5}$$

then $x_0$ is called *fixed point* of the map $f_a$. In the cobweb diagram, fixed points are intersection points of the graph of the function and the straight line $x_{n+1} = x_n$. For the logistic map with $a = 4$, there are only two fixed points: $x_0 = 0$ and $x_0 = 0.75$ (Fig. 6).

Consider now Eq. (3.4) with $n > 1$.

If $x_0$ verifies

$$x_0 = f_a^n(x_0) \tag{3.6}$$

And if

$$f_a^k(x_0) \neq x_0 \quad \text{for } 1 \leq k < n \tag{3.7}$$

then the set $\{x_0, x_1, ..., x_{n-1}\}$, is called the *periodic orbit of period n* of the map. The abbreviated equivalent terms period-$n$ point, period-$n$ orbit; $n$-cycle are often used.

The set $\left\{\dfrac{5+\sqrt{5}}{8}, \dfrac{5-\sqrt{5}}{8}\right\}$ is the period-2 orbit of the logistic map with $a = 4$. However, the period-2 orbit is not the only one periodic orbit for this map (with $a = 4$). In fact there exist an infinity of periodic orbits and an infinity of periods (furthermore several distinct periodic orbits having the same period can coexist).

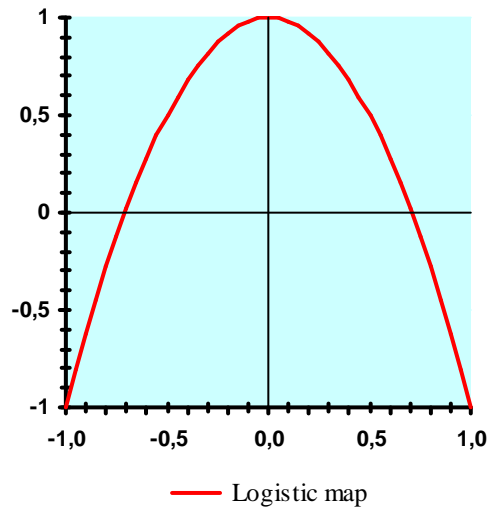

**Fig. 7. Graph of the logistic map** $f(x) = 1 - 2\,x^2$.

It is surprising that a simple quadratic equation can exhibit such complex behaviour. If the logistic equation with $a = 4$ modelled the growth of bugs, then their population would exhibit erratic yearly fluctuations. This map written in the equivalent form in the interval [-1, 1] (Fig. 7)

$$x_{n+1} = 1 - 2x_n^2 \qquad (3.8)$$

was studied by Ulam[4] well before the modern chaos era. He and von Neumann proposed it as a computer random number generator [15].

### 3.2. Probability Distribution, Sharkovskiĭ's Theorem

*(a) Probability Distribution of the Logistic Map with a = 4*

With enough iterations of the logistic map for $a = 4$, the orbit approaches arbitrarily close to every point in the interval $0 \leq x \leq 1$. However, the points visited to the unit interval are not uniformly distributed. In particular, many values in the vicinity of $x_n = 0$ map into values of $x_{n+1}$ close to one, which in turn map into values of $x_{n+2}$ close to zero. Hence the *probability distribution function* (also called the *probability density distribution*, the *invariant measure* or the *natural measure*) $P(x)$ has peaks at $x = 0$ and $x = 1$. $P(x)$ is the probability that a point $x$ is within $dx$ of $x$ and is normalized so that the area under the $P(x)$ curve is

$$\int_0^1 P(x)dx = 1 \tag{3.9}$$

The point has the unit probability of being somewhere in the interval $0 \leq x \leq 1$. $P(x)$ can be explicitly calculated [11, 12]. The result is

$$P(x) = \frac{1}{\pi\sqrt{x(1-x)}} \tag{3.10}$$

as shown in Fig. 8.

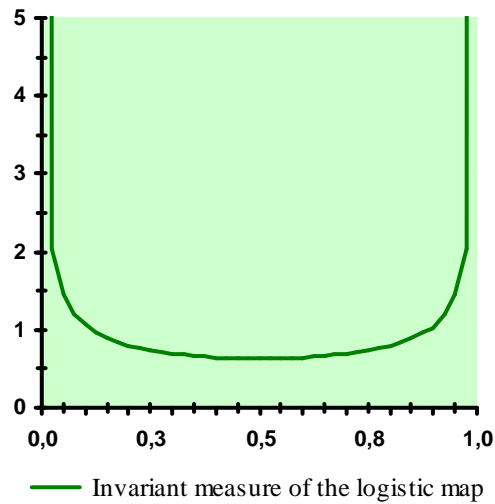

— Invariant measure of the logistic map

**Fig. 8. Graph of the invariant measure** $P(x) = \frac{1}{\pi\sqrt{x(1-x)}}$ **of the logistic map** $f(x) = 4\,x\,(1-x)$.

Practically, if we use a finite number $N$ of iterated points of the logistic map, only an approached probability distribution $P_N(x)$ is obtained (Fig. 9). As the number of iterates increases, $P_N(x)$ is converging toward $P(x)$. In next section we will consider the difference between $P_N(x)$ and $P(x)$ as a benchmark of the model of coupled logistic maps we introduce.

*(b) Sharkovskiĭ's Theorem*

A continuous map of the unit interval [0, 1] may have one fixed point and no other periodic orbits (for example, $f(x) = x/2$). There may be fixed points, period-two orbits, and no other periodic orbits (for example, $f(x) = 1 - x$). More complicated situations may occur. In 1975, T.-Y. Li and J. A. Yorke proved that "period three implies chaos" [16], that means that if a three-periodic orbit exists for a continuous map of the interval, then for every positive natural number $k$ there exists at

---

[4] S. T. Ulam (1909 - 1984) a Polish mathematician emigrated to the United States was instrumental in developing the

least one *k*-periodic orbit.

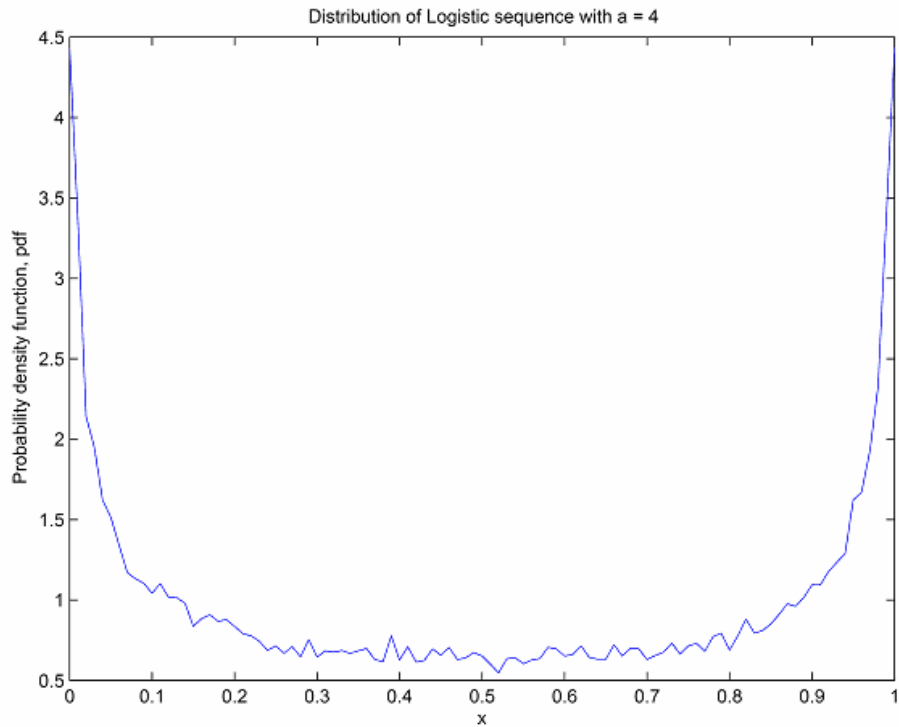

**Fig. 9**. **Distribution of logistic sequence in unstable state ($a = 4$) (from [4]).**

In fact, a more general theorem concerning the existence of periodic orbits of continuous maps of the real line into itself was proved thirteen years before in 1962 by Sharkovskiĭ [17], an Ukrainian mathematician[5]. However, this theorem was not easily accessible in the western literature and remained unknown to the English speaking math community until 1976.

This theorem may be stated as follows:

Consider the set of natural numbers in which the following relation has been introduced: $n_1$ precedes $n_2$ ($n_1 \prec n_2$) if *for every continuous map of the line into itself the existence of a cycle of order $n_1$ implies the existence of a cycle of order $n_2$*. This relation is clearly reflexive and transitive and, consequently, the set of natural numbers with this relation is a quasi-ordered set[6].

---

hydrogen bomb and the *Monte Carlo method* for evaluating complicated integrals.

[5] Alexander Nikolaevitch Sharkovskiĭ (1936- ) is a Ukrainian mathematician who received his PhD from the Institute of Mathematics of the National Academy of Sciences of Ukraine in 1961, of which he is corresponding member since 1978. It has been established that the Sharkovskiĭ's order or its modifications occur in many other systems (among which are random and multivalent systems, many-dimensional and even infinite-dimensional systems), and not only for cycles but for more complicated structures as well.

[6] In mathematics, a well-order (or well-ordering) on a set S is a total order on S with the property that every non-empty subset of S has a least element in this ordering. The set S together with the well-order is then called a *well-ordered* set. Roughly speaking, a well-ordered set is ordered in such a way that its elements can be considered one at a time, in order, and any time you have not examined all of the elements, there is always a unique next element to consider (Wikipedia encyclopedia).

Sharkovskiĭ 's theorem: This relation turns the set of natural numbers into *n* ordered set, which is ordered in the following way:

$$3 \prec 5 \prec 7 \prec 9 \prec 11 \prec \cdots \prec 3 \cdot 2 \prec 5 \cdot 2 \prec \cdots \prec 3 \cdot 2^2 \prec 5 \cdot 2^2 \prec \cdots \prec 2^3 \prec 2^2 \prec 2^1 \prec 1 \quad (3.11)$$

We can do some remarks:

**Remark 1.** If *f* has a periodic point whose period is not a power of two, then *f* must have infinitely many periodic points.

**Remark 2.** Sharkovskiĭ's theorem is sharp. In the logistic map when $1 + \sqrt{8} < a < 3.8415....$ *f* has a 3- cycle. It must also have cycle of any period. The value $1 + \sqrt{8} = 3.828427...$ being explicitly calculated, because it is the critical value *a*\* for the period 3 tangent bifurcation.

**Remark 3**. Sharkovskiĭ's theorem tells us nothing about the stability of the orbits, or the parameter values for which they may be observed.

However, in some cases, depending on the nature of the map considered, it is possible to obtain some results on the nature of the orbit. For example, in the logistic map, in spite of the existence of infinitely many periodic orbits, at most one is attracting, because the theorem of Fatou says that, "every attracting cycle for a polynomial or a rational function attracts at least one critical point."

For the logistic map, it has only one critical point *x* = 1/2, were *f'*(1/2) = 0. Thus, for every value of *a*, *f* has at most one attracting cycle.

### 3.3. Orbits of Discrete Approximations to Chaotic Maps

It should be impossible to study the behaviour of chaotic maps without the help of computers. Of course, some theoretical results as Sharkovskiĭ's theorem or the exact value of the probability distribution of the logistic map for *a* = 4 and some other important results are obtained without this help [18, 19], however it is noticeable that the most known examples of chaotic attractor were found using numerical simulation (i.e. the Lorenz attractor [20] in 1963 and the Hénon attractor [21] in 1976). The use of computers allows us to explore the behaviour of iterates for large values of *n*. However, strong limitations of this use appear when exact theorems are needed. This is due to the nature of the numbers handled by computers. We recall some features of such numbers below.

*(a) Floating Point Numbers*

There are several ways to represent real numbers on computers [22, 23]. Fixed point places a radix point somewhere in the middle of the digits, and is equivalent to using integers that represent portions of some unit. For example, one might represent 1/100 ths of a unit; if one has four decimal digits, one could represent 27.46, or 00.07. Another approach is to use rationals, and represent every number as the ratio of two integers e.g. 27/56.

*Floating-point representation*, the most common solution, basically represents reals in scientific notation. Scientific notation represents numbers as a base number and an exponent. For example, 123.456 could be represented as $1.23456 \times 10^2$. In hexadecimal, the number 456.ABC might be represented as $4.56ABC \times 10^2$ or more explicitly $4.56ABC_{hexa} \times 16^2_{decimal}$.

Floating-point solves a number of representation problems. Fixed-point has a fixed window of representation, which limits it from representing very large or very small numbers. Also, fixed-point is prone to a loss of precision when two large numbers are divided.

Floating-point, on the other hand, employs a sort of "sliding window" of precision appropriate to the scale of the number. This allows it to represent numbers from 2,000,000,000,000 to 0.0000000000000002 with ease.

IEEE Standard 754 floating point is the most common representation today for real numbers on computers, including Intel-based PC's, Macintoshes, and most Unix and Linux platforms.

IEEE floating point numbers have three basic components, viz. the sign, the exponent, and the mantissa. The mantissa is composed of the *fraction* and an implicit leading digit (explained below). The exponent base (**2**) is implicit and need not be stored.

Table 1 shows the layout for single (32-bit) and double (64-bit) precision floating-point values. The number of bits for each field are shown (bit ranges are in square brackets).

**Table 1. Layout for single (32-bit) and double (64-bit) precision floating-point values**

|  | Sign | Exponent | Fraction | Bias |
|---|---|---|---|---|
| **Single Precision** | 1 [31] | 8 [30-23] | 23 [22-00] | 127 |
| **Double Precision** | 1 [63] | 11 [62-52] | 52 [51-00] | 1023 |

*(i) The Sign Bit*

The sign bit is as simple as it gets. 0 denotes a positive number, 1 denotes a negative number. Flipping the value of this bit flips the sign of the number.

*(ii) The Exponent*

The exponent field needs to represent both positive and negative exponents. To do this, a bias is added to the actual exponent in order to get the stored exponent. For IEEE single-precision floats, this value is 127. Thus, an exponent of zero means that 127 is stored in the exponent field. A stored value of 200 indicates an exponent of (200-127), or 73. For double precision, the exponent field is 11 bits, and has a bias of 1023.

*(iii) The Mantissa*

The *mantissa*, also known as the *significand*, represents the precision bits of the number. It is composed of an implicit leading bit and the fraction bits.

To find out the value of the implicit leading bit, consider that any number can be expressed in scientific notation in many different ways. For example, the number seven can be represented as any of these:

$$7.00 \times 10^0$$
$$0.007 \times 10^3$$
$$70000 \times 10^{-4}$$

In order to maximize the quantity of representable numbers, floating-point numbers are typically stored in *normalized* form. This basically puts the radix point after the first non-zero digit. In normalized form, seven is represented as $7.0 \times 10^0$.

A nice little optimization is available to us in base two, since the only possible non-zero digit is 1. Thus, we can just assume a leading digit of 1, and do not need to represent it explicitly. As a result, the mantissa has effectively 24 bits of resolution, by way of 23 fraction bits.

*(iv) Ranges of Floating Point Numbers*

Let's consider single-precision floats for a second. Note that we are taking essentially a 32-bit number and re-jiggering the fields to cover a much broader range. Something has to give, and it is precision. For example, regular 32-bit integers, with all precision centered around zero, can precisely store integers with 32-bits of resolution. Single-precision floating-point, on the other hand, is unable to match this resolution with its 24 bits. It does, however, approximate this value by effectively truncating from the lower end.

For example

11110000 11001100 10101010 0000**1111**  // 32-bit integer

= +1.1110000 11001100 10101010 x $2^{31}$   // Single-Precision Float

=  11110000 11001100 10101010 0000**0000**  // Corresponding Value

This approximates the 32-bit value, but does not yield an exact representation. The difference is bold typed. On the other hand, besides the ability to represent fractional components (which integers lack completely), the floating-point value can represent numbers around $2^{127}$, compared to 32-bit integers maximum value around $2^{32}$.

The effective range (excluding infinite values) of IEEE floating-point numbers is shown in Table 2.

Since the sign of floating point numbers is given by a special leading bit, the range for negative numbers is given by the negation of the above values.

**Table 2. Effective range of positive floating point numbers**

|  | Binary | Decimal |
|---|---|---|
| **Single** | $\pm (2-2^{-23}) \times 2^{127}$ | $\sim \pm 10^{38.53}$ |
| **Double** | $\pm (2-2^{-52}) \times 2^{1023}$ | $\sim \pm 10^{308.25}$ |

*(b) Computational Experiments*

When a dynamical system is realized on a computer using floating point numbers described as above, the computation is of a discretization, where finite machine arithmetic replaces continuum state space. For chaotic dynamical systems, the discretizations often have collapsing effects to a fixed point or to short cycles.

*(i) Approximated Logistic Map*

As an example of such collapsing effects, Lanford III [24], presents the results of a sampling study in double precision of a discretization of the logistic map (3.8) which has excellent ergodic properties. The precise discretization studied is obtained by first exploiting evenness to fold the interval [-1,0] to [0,1], i.e., replacing (3.8) by

$$x_{n+1} = \left| 1 - 2x_n^2 \right| \qquad (3.12)$$

On [0, 1] it is not difficult to see that the folded map has the same periods as the original one. The working interval is then shifted from [0, 1] to [1, 2] by translation, in order to avoid perturbation in numerical experiment, by the special value 0. Then the translated folded map is programmed in straightforward way. Out of 1 000 randomly chosen initial points:

- 890, i.e., the overwhelming majority, converged to the fixed point corresponding to the fixed point -1 in the original representation (3.8).
- 108 converged to a cycle of period 3,490,273.
- the remaining 2 converged to a cycle of period 1,107,319.

Thus, in this case at least, the very long-term behaviour of numerical orbits is, for a substantial fraction of initial points, in flagrant disagreement with the true behaviour of typical orbits of the original smooth logistic map.

*(ii) Another Approximated Maps*

Lanford III [24] carefully studies the numerical approximation of the map

$$x_{n+1} = 2x_n + 0.5x_n(1-x_n) \quad \text{(mod 1)} \qquad 0 \leq x \leq 1 \qquad (3.13)$$

It is perhaps better to think this map (see Fig. 10) as acting on the unit interval with endpoints identified, i.e., on the circle. Note that $f'(x) \geq 1.5$ everywhere, so $f$ is strictly expanding in a particularly clean and simple sense. As a consequence of expansivity, this mapping has about the best imaginable ergodic properties:

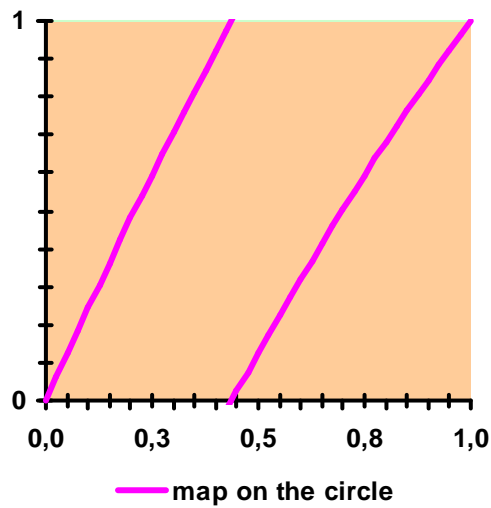

**Fig. 10.** **Graph of the map** $f(x) = 2x + 0.5\,x\,(1-x)$.

- it admits a unique invariant measure µ equivalent to Lebesgue measure,
- the abstract dynamical system $(f, \mu)$ is ergodic and in fact isomorphic to a Bernoulli shift,
- a central limit theorem holds, …

One consequence of ergodicity of $f$ relative to µ is that, for Lebesgue almost initial points $x$ in the unit interval, the corresponding orbit $f^{(n)}(x)$ is asymptotically distributed over the unit interval according to µ.

In numerical experiments performed by Lanford III, the computer working with fixed finite precision is able to represent finitely many points in the interval in question. It is probably good, for purposes of orientation, to think of the case where the representable points are uniformly spaced in the interval. The true smooth map is then *approximated* by a discretized map, sending the finite set of representable points in the interval to itself.

Describing the discretized mapping exactly is usually complicated, but it is *roughly* the mapping obtained by applying the exact smooth mapping to each of the discrete representable points and "rounding" the result to the nearest representable point[7].

He has done two kinds of experiments: first, using uniformly spaced points in the interval with several order of discretization (ranging from $2^{22}$ to $2^{25}$), and second, using double precision floating points.

In each experiment the questions addressed are:
- how many periodic cycles are there and what are their periods ?

---

[7] The reason why this simple description is not quite realistic is that, in practice, intermediate quantities, and not just the final result, undergo rounding.

- how large are their respective basins of attraction, i.e., for each periodic cycle, how many initial points give orbits with eventually land on the cycle in question ?

For relatively coarse discretizations the orbit structure is determined completely, i.e., all the periodic cycles and the exact sizes of their basins of attraction are found. Some representative results are given in Tables 3 to 7. In theses tables, **N** denotes the order of the discretization, i.e., the representable points are the numbers, $j/N$, with $0 \leq j < N$.

**Table 3. Coexisting periodic orbits of mapping (3.13) for the discretization $N = 2^{22}$.**

| $N = 2^{22} = 4{,}194{,}304$ | | |
| --- | --- | --- |
| 13 cycles | | |
| Period | Basin size | Relative size |
| 3,864 | 2,523,929 | 60.18 % |
| 1,337 | 538,712 | 12.84 % |
| 718 | 513,839 | 12.26 % |
| 295 | 238,486 | 5.69 % |
| 130 | 203,587 | 4.86 % |
| 1,338 | 152,942 | 3.65 % |
| 297 | 12,359 | 0.29 % |
| 169 | 5,056 | 0.12 % |
| 97 | 3,012 | 0.07 % |
| 17 | 2,346 | 0.06 % |
| 6 | 21 | $5.10^{-4}$ % |
| 1 | 8 | $2.10^{-4}$ % |
| 7 | 7 | 5.10 |

**Table 4. Coexisting periodic orbits of mapping (3.13) for the discretization $N = 2^{23}$.**

| $N = 2^{23} = 8{,}388{,}608$ | | |
| --- | --- | --- |
| 7 cycles | | |
| Period | Basin size | Relative size |
| 4,898 | 5,441,432 | 64.85 % |
| 1,746 | 2,946,734 | 35.13 % |
| 13 | 205 | $24.10^{-4}$ % |
| 6 | 132 | $16.10^{-4}$ % |
| 30 | 96 | $11.10^{-4}$ % |
| 4 | 8 | $< 1.10^{-6}$ % |
| 1 | 1 | $< 1.10^{-6}$ % |

**Table 5. Coexisting periodic orbits of mapping (3.13) for the discretization $N = 2^{24}$.**

| $N = 2^{24} = 16{,}777{,}216$ | | |
| --- | --- | --- |
| 2 cycles | | |
| Period | Basin size | Relative size |
| 5,300 | 16,777,214 | 100 % |
| 1 | 2 | $< 1.10^{-6}$ % |

**Table 6. Coexisting periodic orbits of mapping (3.13) for the discretization $N = 2^{24}-1$.**

| $N = 2^{24} -1= 16,777,215$ 10 cycles | | |
|---|---|---|
| Period | Basin size | Relative size |
| 3,081 | 7,502,907 | 44.72 % |
| 699 | 3,047,369 | 18.16 % |
| 3,469 | 2,905,844 | 17.32 % |
| 1,012 | 2,774,926 | 16.54 % |
| 563 | 290,733 | 11.73 % |
| 2,159 | 221,294 | 1.32 % |
| 138 | 21,610 | 0.13 % |
| 421 | 12,477 | 0.07 % |
| 9 | 54 | $< 1.10^{-3}$ % |
| 1 | 1 | $< 1.10^{-7}$ % |

**Table 7. Coexisting periodic orbits of mapping (3.14) for the discretization $N = 2^{25}$.**

| $N = 2^{25} = 33,554,432$ 8 cycles | | |
|---|---|---|
| Period | Basin size | Relative size |
| 4,094 | 32,114,650 | 95.71 % |
| 621 | 918,519 | 2.74 % |
| 283 | 516,985 | 1.54 % |
| 126 | 2,937 | $< 0.01$ % |
| 6 | 887 | $< 0.01$ % |
| 55 | 433 | $< 0.01$ % |
| 4 | 20 | $< 1.10^{-6}$ % |
| 1 | 1 | $< 1.10^{-6}$ % |

For the experiments using ordinary (IEEE-754) double precision, so that the working interval contains of the order of $10^{16}$ representable points, 1000 initial points at random are used in order to sample the orbit structure, determining the cycles to which they converge.

The computations were accomplished by shifting the working interval from [0, 1] to [1, 2] by translation, for the same previous reason as for the logistic map, i.e. , the map actually iterated was

$$x_{n+1} = 2x_n + 0.5x_n(x_n - 1)(2 - x_n) \quad (\text{mod } 1) \qquad 1 \leq x \leq 2 \qquad (3.14)$$

Some results are displayed in Table 8.

Many more examples could be given, but those given may serve to illustrate the intriguing character of the results; the outcomes proves to be extremely sensitive to the details of the experiment, but the results all have a similar flavour: a relatively small number of cycles attract near all orbits, and the lengths of these significant cycles are much larger than one but much smaller than the number of representable points. Lanford III wrote that there are here regularities which ought to be understood.

**Table 8. Coexisting periodic orbits of mapping (3.13) for the discretization N = $2^{25}$.**

| Double precision (sampling) 7 cycles found | | |
|---|---|---|
| Period | Basin size | Relative size |
| 27,627,856 | 517 | 51.7 % |
| 88,201,822 | 320 | 32.0 % |
| 4,206,988 | 147 | 14.7 % |
| 4,837,566 | 17 | 1.7 % |
| 802,279 | 8 | 0.8 % |
| 6,945,337 | 6 | 0.6 % |
| 2,808;977 | 1 | 0.1 % |

Diamond and Pokrovskii [25] suggest that statistical properties of the phenomenon of computational collapse of discretized chaotic mapping can be modelled by random mappings with an absorbing centre. The model gives results which are very much in line with computational experiments and there appears to be a type of universality summarised by an Arcsine law. The effects are discussed with special reference to the family of mappings

$$x_{n+1} = 1 - |1 - 2x_n|^l \quad 0 \leq x \leq 1 \quad 1 < l \leq 2 \quad (3.15)$$

Computer experiments show close agreement with prediction of the model.

### 3.4. Mega and Gigaperiodic Orbits

It should be noticed that several orbits observed in both approximated maps (3.12) (logistic) or (3.13-3.14) (circle map) have lengths of order $10^6$ or $10^7$. We propose to give the name "Megaperiodic" to these orbits in reference to the prefix Mega which is for $10^6$.

More precisely we will call *megaperiodic* orbits, those whose length of the periods belong to the interval of natural numbers $[10^6, 10^9[$. In the same way we will call *gigaperiodic* orbits, those whose length of the periods belong to the interval $[10^9, 10^{12}[$.

In Table 9, we define precisely mega and giga periodic orbits and also, tera, peta, exa, zetta, yottaperiodic orbits with respect to the denomination of the international systems of prefixes of physical units.

**Remark.** When computers were running slower than now in 1976, Michel Hénon (mathematician and astronomer at the observatory of Nice, France) pointed out numerically the first example of strange attractor of a map of the plane [21]. Using an IBM 7040 computer with 16-digit accuracy he plotted some of the first $5.10^6$ iterates of both the initial points (0, 0) and (0.63135448, 0.18940634) in order to display graphically the fractal structure of his attractor.

**Table 9.** Definition of mega, giga, tera , … periodic orbits with respect to the length of their periods.

| Name | Symbol | Length of the period |
|---|---|---|
| Megaperiodic orbit | **M**-periodic | $10^6 \leq$ period $< 10^9$ |
| Gigaperiodic orbit | **G**-periodic | $10^9 \leq$ period $< 10^{12}$ |
| Teraperiodic orbit | **T**-periodic | $10^{12} \leq$ period $< 10^{15}$ |
| Petaperiodic orbit | **P**-periodic | $10^{15} \leq$ period $< 10^{18}$ |
| Exaperiodic orbit | **E**-periodic | $10^{18} \leq$ period $< 10^{21}$ |
| Zettaperiodic orbit | **Z**-periodic | $10^{21} \leq$ period $< 10^{24}$ |
| Yottaperiodic orbit | **Y**-periodic | $10^{24} \leq$ period $< 10^{27}$ |

The discovery of a strange attractor was only possible because, thirty years ago, nobody was interested in searching long period (**G**-period). Then it was unusual to consider periodic orbits with periods longer than 10 or 20. As said in [21] " … *a limit cycle intersecting the surface of section **p** times. The value of **p** increases through successive "bifurcations" as **a** increases, and appears to tend to infinity as **a** approaches a critical value of the order 1.06 …*".

In fact using now a desktop computer with a 3 Gigahertz processor, allows us to compute circa 20,000,000 iterates of the Hénon map per second. Hence it is very easy to look for mega or gigaperiodic orbits, and the "numerical strange attractor" found by him, is no more than a periodic orbit with a particular geometrical structure.

### 3.5. Chaotic Maps of the Plane

Until now we have emphasized one dimensional dynamical systems. However dynamical systems in dimension two or more are much more interesting, because they produce more complex behaviour and also periodic orbits with larger period. There are several chaotic maps of the plane including baker's map [26], Belykh map [27], Cohen map [28], Duffing map, Gaussian map, Hénon map [21], Kaplan-Yorke map [29], Lozi map [30], sinusoidal map, standard map [31], Zaslavskii map [32], …and also a huge number of publications on these chaotic maps and their mathematical properties. Afraimovich et al. [33] proposed to divide the chaotic attractors of these maps into robust hyperbolic attractors, almost hyperbolic (quasihyperbolic) attractors and so-called quasiattractors. Such a classification of attractors is a natural consequence of the rigorous mathematical analysis of the structure and properties of dynamical chaos. However, since from the experimental point of view it has not been accepted as significant, Anishchenko et al. [34] studied quasihyperbolic attractors and quasiattractors in two-dimensional invertible maps in order to obtain some characteristic properties which allow one to diagnose exactly their difference.

As pointed out in Section 1, the difference between chaos and random has been recently carefully studied and for several applications chaos has been found more useful than random in some cases [3, 4]. Both Hénon and Lozi map are now often used in order to generate chaotic sequences. In [35] the time series of these maps are used to test a minimax method for learning functional networks. Functional networks are a soundness and efficient generalization of neural networks. In [36], Lozi sequences are used to show the good performance of feedforward neural networks and functional networks to model nonlinear chaotic systems. Park et al. [37] test a new generalized predictive control method based on an ARMAX model with these chaotic sequences. Caponetto et al. [38] introduce these chaotic sequences instead of random ones in order to improve the performance of evolutionary algorithms. Borges et al. [39] show that the sensitivity to initial conditions of both maps (in fact a smooth version of the Lozi map) is the same. Anishchenko et al. [40] point out that the Lozi attractor and the Lorenz attractor behave in a similar manner, moreover their characteristics measured in numerical experiments, are robust relative to small perturbations.

### 3.5.1. The Hénon map

The Hénon map is an invertible mapping of a two-dimensional plane into itself

$$\begin{cases} x_{n+1} = y_n + 1 - a x_n^2 \\ y_{n+1} = b x_n \end{cases} \qquad (3.16)$$

Equivalently, the Hénon map can be defined by the 2-step recurrence relation, also called 2D *delayed map* of the form $x_{n+1} = F(x_n, x_{n-1})$

$$x_{n+1} = 1 - a x_n^2 + b x_{n-1} \qquad (3.17)$$

The Hénon map is one of the simplest models of a Poincaré map of a three-dimensional invertible flow. For $a = 1.4$ and $b = 0.3$ there is a strange attractor (Fig. 11).

*(a) Two G-periodic Orbits*

In comparison with the two **M-**periodic orbits displayed on Table 8 for the map (3.14), Hénon map has **G-**periodic orbits. On a Dell computer with a Pentium IV microprocessor running at the frequency of 1.5 Gigahertz, using a Borland C compiler and computing with ordinary (IEEE-754) double precision numbers, we can find for $a = 1.4$ and $b = 0.3$ an attracting period of length **3,800,716,788,** i.e., forty times longer than the longest period of the one-dimensional map.

This periodic orbit (we call it here Orbit 1) is numerically slowly attracting: starting with the initial value

$$(x_0, y_0)_1 = (-0.35766, 0.14722356)$$

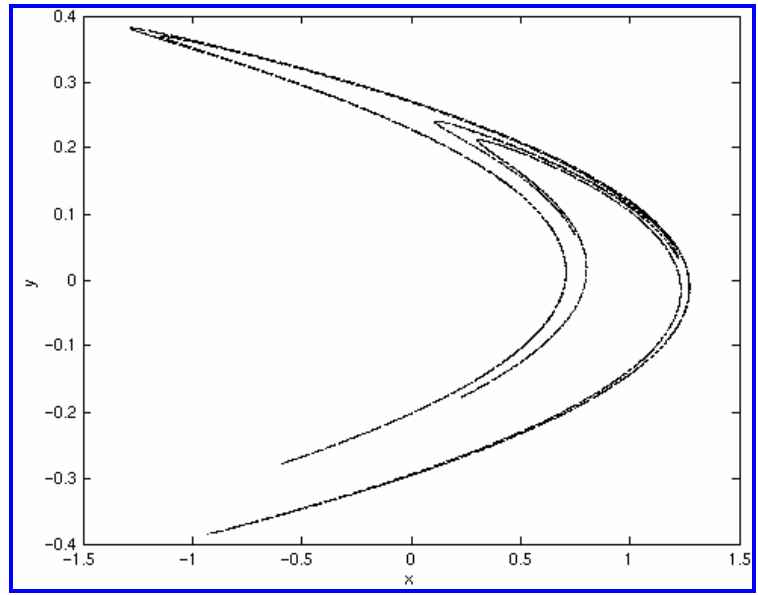

**Fig. 11.** **Hénon map, plot of the first 5,000 iterations using** $a = 1.7$ **and** $b = 0.5$.

We obtain

$(x_{11,574,730,767}, y_{11,574,730,767})_1 = (1.272973613509556 62, -0.0115735710153616837)$

and

$(x_{15,375,447,555}, y_{15,375,447,555})_1 = (1.27297361350955662, -0.0115735710153616837)$

Subtracting, the length of the period is obtained as

$$15{,}375{,}447{,}555 - 11{,}574{,}730{,}767 = 3{,}800{,}716{,}788.$$

However, as for the map (3.14) studied by Lanford, this periodic orbit is not unique: starting with the initial value

$$(x_0, y_0)_2 = (0.4725166, 0.25112222222356)$$

We obtain

$(x_{12,935,492,515}, y_{12,935,492,515})_2 = (x_{13,246,439,123}, y_{13,246,439,123})_2$

$= (1.27297361350865113, -0.0115734779561870744)$

which is a **M**-periodic orbit of period 310,946,608 (Orbit 2).

This orbit can be reached more rapidly starting from the other initial value

$$(x_0, y_0) = (0.881877775591, 0.0000322222356)$$

then $(x_{4,459,790,707}, y_{4,459,790,707}) = (1.27297361350865113, -0.0115734779561870744)$.

It is possible that some others periodic orbits coexist with both orbit 1 and orbit 2.

*(b) Sensitivity to Initial Condition*

The comparison between orbit 1 and orbit 2 gives a perfect idea of the sensitive dependence on initial conditions of chaotic attractors:

Orbit 1 passes through the point

$$(\mathbf{1.272973613 50}955662, \mathbf{-0.0115 73}5710153616837)$$

and orbit 2 passes through the point

$$(\mathbf{1.272973613 50}865113, \mathbf{-0.0115 73}4779561870744)$$

The same digits of these points are bold printed, they are very near.

**Remark.** For the Hénon map, in order to found fixed points, we have to solve the quadratic equation

$$ax^2 + (1-b)x - 1 = 0 \tag{3.18}$$

To obtain periodic solution we have to solve polynomial equation with degree larger than 2, which never be done explicitly, except for 2-cycle.

### 3.5.2. The Lozi map

The Lozi map is a linearized version of the Hénon map, built in order to simplify the computations, mainly because it is possible to compute explicitly any periodic orbits solving a linear system.

$$\begin{cases} x_{n+1} = y_n + 1 - a|x_n| \\ y_{n+1} = bx_n \end{cases} \tag{3.19}$$

or equivalently

$$x_{n+1} = 1 - a|x_n| + bx_{n-1} \tag{3.20}$$

For $a = 1.7$ and $b = 0.5$ (Fig. 12) there is a strange attractor. The particularity of this strange attractor is that it has been rigorously proved by Misiurewicz in 1980 [41].

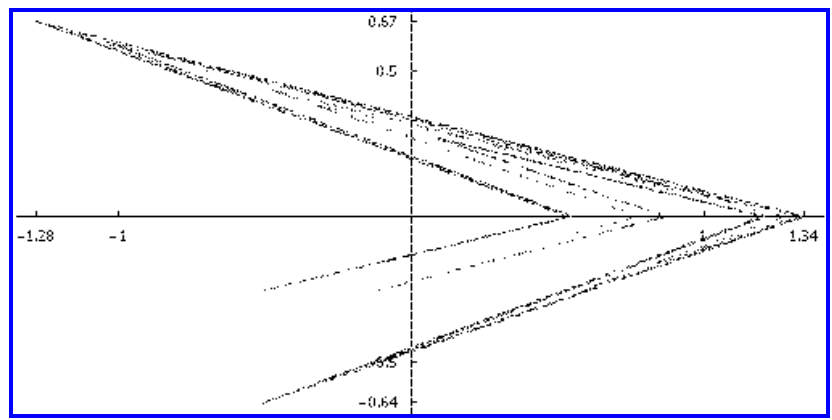

**Fig. 12.**   Lozi map, plot of the first 5,000 iterations using $a = 1.7$ and $b = 0.5$.

*(a) One G-periodic Orbit*

In the same conditions of computation, running the computation during 19 hours, we can find a **G**-periodic attracting orbit of period **436,170,188,959.**

Starting with

$$(x_0, y_0) = (0.88187777591, 0.0000322222356)$$

We obtain

$$(x_{686,295,403,186}, y_{686,295,403,186}) = (x_{250,125,214,227}, y_{250,125,214,227})$$
$$= (1.34348444450739479, -2.07670041497986548 \cdot 10^{-7})$$

There is a transient regime before the orbit is reached.

It seems that there is no periodic orbit with a smaller length, this could be due to the quasihyperbolic nature of the attractor. Hence, this attractor is very efficient, in order to generate chaotic numbers without repetition for standard simulation. However they are not *equally distributed* on the plane (Fig. 12) (or at least in the domain of attraction of the attractor in the plane).

To improve the generation of chaotic numbers we will consider in the next section how to generate chaotic numbers with uniform repartition on a given interval, or on a given square of the plane or more generally in a given hypercube of $\mathbf{R}^n$.

## 4. Very Weakly Coupled Symmetric Tent Maps

We introduce a new point of view in the study of coupled maps. For instance, the system of two coupled maps of the interval [42] is

$$\begin{cases} x_{n+1} = (1-\varepsilon)f(x_n) + \varepsilon f(y_n) \\ y_{n+1} = \varepsilon f(x_n) + (1-\varepsilon)f(y_n) \end{cases} \tag{4.1}$$

The coupling constant $\varepsilon$ varies from 0 to 1. When $\varepsilon = 0$, the maps are decoupled, when $\varepsilon = 1$ they are fully cross coupled. Generally, numerical computations do not consider very very small values of $\varepsilon$ (as small as $10^{-7}$ for floating numbers or $10^{-14}$ for double precision numbers), because it seems that the maps are nearly decoupled with those values. Hence, no special effect of the coupling is expected.

We will see in this section that it is not the case and that very very small coupling constant allows us the construction of very long periodic orbits, leading to sterling chaotic generators.

### 4.1. The Symmetric Tent Map

In this paper we only consider the logistic map (3.2) with $a = 2$ and the symmetric tent map as component of these new chaotic generators (Fig. 13). The symmetric tent map is defined by

$$f_a(x) = 1 - a|x| \tag{4.2}$$

We consider also this map with the value $a = 2$.

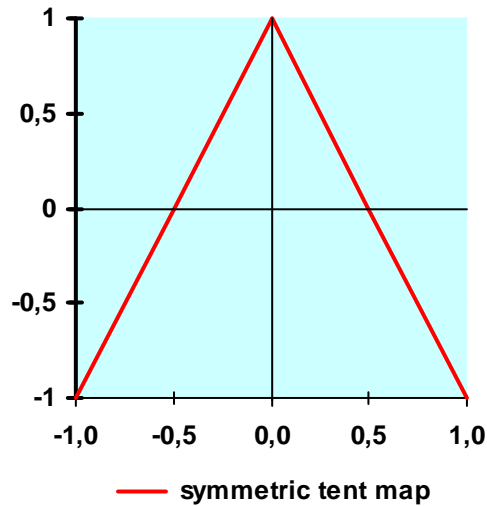

**Fig. 13.** **Graph of the symmetric tent map** $f(x) = 1 - 2|x|$.

The associated dynamical system is defined by the equation on the interval [-1, 1]

$$x_{n+1} = 1 - a|x_n| \qquad (4.3)$$

Despite its simple shape, it has several interesting properties. First, it has chaotic orbits. Because of its simple shape, it is easy to find those orbits explicitly.

Also because of its simple shape, the tent map's shape under iteration is very well understood. The invariant measure is the Lebesgue measure. Finally, and perhaps most importantly, the tent map is conjugate to the logistic map, which in turn, is conjugate to the Hénon map for small values of *b*.

However, the symmetric tent map is dramatically numerically instable: Sharkovskiĭ's theorem applies for it. There is a period three orbit which implies that there is infinity of periodic orbits. Nevertheless the orbit of almost every point of the interval [-1, 1] of the discretized tent map converges to the (unstable) fixed point -1, and there is no numerical attracting periodic orbit.

This is why the tent map is never used to generate numerically chaotic numbers.

The behaviour of iterates with respect to chaos is worse than the behaviour of iterates of the approximated logistic map.

**4.2. Two Coupled Maps**

Our aim now is to increase the length of the periods of the map of the plane we introduce, coupling two identical maps of the interval. We generalize equation (4.1) allowing us to use a two dimensional coupling constant $\varepsilon = (\varepsilon_1, \varepsilon_2)$:

$$\begin{cases} x_{n+1} = (1-\varepsilon_1)f(x_n) + \varepsilon_1 f(y_n) \\ y_{n+1} = \varepsilon_2 f(x_n) + (1-\varepsilon_2)f(y_n) \end{cases} \quad (4.4)$$

However, in this paper we fix constant the ratio between $\varepsilon_1$ and $\varepsilon_2$ equal to 2. Others ratios can also lead to good results.

$$\varepsilon_2 = 2\varepsilon_1 \quad (4.5)$$

The dynamical system (4.4) can be described more generally by

$$\begin{pmatrix} x_{n+1} \\ y_{n+1} \end{pmatrix} = F\begin{pmatrix} x_n \\ y_n \end{pmatrix} = \begin{pmatrix} (1-\varepsilon_1) & \varepsilon_1 \\ \varepsilon_2 & (1-\varepsilon_2) \end{pmatrix} \cdot \begin{pmatrix} f(x_n) \\ f(y_n) \end{pmatrix} \quad (4.6)$$

Or

$$X_{n+1} = F(X_n) = A \cdot \underline{f}(X_n) \quad (4.7)$$

With

$$X = \begin{pmatrix} x \\ y \end{pmatrix}, \quad \underline{f}(X) = \begin{pmatrix} f(x) \\ f(y) \end{pmatrix} \quad \text{and} \quad A = \begin{pmatrix} (1-\varepsilon_1) & \varepsilon_1 \\ \varepsilon_2 & (1-\varepsilon_2) \end{pmatrix}$$

$F$ is a map of the square [-1, 1] x [-1, 1] into itself.

### 4.3. Three and *p*-coupled Symmetric Maps

In order to improve the length of the period and the convergence of the invariant measure to a given measure, we will have to consider the coupling of three or more maps of the interval with the following hypothesis:

$$\begin{cases} x_{n+1} = (1-2\varepsilon_1)f(x_n) + \varepsilon_1 f(y_n) + \varepsilon_1 f(z_n) \\ y_{n+1} = \varepsilon_2 f(x_n) + (1-2\varepsilon_2)f(y_n) + \varepsilon_2 f(z_n) \\ z_{n+1} = \varepsilon_3 f(x_n) + \varepsilon_3 f(y_n) + (1-2\varepsilon_3)f(z_n) \end{cases} \quad (4.8)$$

and

$$\varepsilon_2 = 2\varepsilon_1, \; \varepsilon_3 = 3\varepsilon_1 \quad (4.9)$$

More generally the model of *p* coupled maps we introduce takes the form

$$\begin{pmatrix} x_{n+1}^1 \\ x_{n+1}^2 \\ \vdots \\ \vdots \\ x_{n+1}^p \end{pmatrix} = F\begin{pmatrix} x_n^1 \\ x_n^2 \\ \vdots \\ \vdots \\ x_n^p \end{pmatrix} = \begin{pmatrix} 1-(p-1)\varepsilon_1 & \varepsilon_1 & \cdots & \varepsilon_1 & \varepsilon_1 \\ \varepsilon_2 & 1-(p-1)\varepsilon_2 & \cdots & \varepsilon_2 & \varepsilon_2 \\ \vdots & & \ddots & \vdots & \vdots \\ \vdots & & & \ddots & \vdots \\ \varepsilon_p & \cdots & \cdots & \varepsilon_p & 1-(p-1)\varepsilon_p \end{pmatrix} \times \begin{pmatrix} f(x_n^1) \\ f(x_n^2) \\ \vdots \\ \vdots \\ f(x_n^p) \end{pmatrix} \quad (4.10)$$

with

$$\varepsilon_i = i\,\varepsilon_1 \qquad i = 1, \ldots, p \quad (4.11)$$

or

$$X_{n+1} = F(X_n) = A \cdot \underline{f}(X_n)$$

with
$$X = \begin{pmatrix} x^1 \\ \vdots \\ \vdots \\ x^p \end{pmatrix}, \quad \underline{f}(X) = \begin{pmatrix} f(x^1) \\ \vdots \\ \vdots \\ f(x^p) \end{pmatrix}$$

and
$$A = \begin{pmatrix} 1-(p-1)\varepsilon_1 & \varepsilon_1 & \cdots & \varepsilon_1 & \varepsilon_1 \\ \varepsilon_2 & 1-(p-1)\varepsilon_2 & \cdots & \varepsilon_2 & \varepsilon_2 \\ \vdots & & \ddots & \vdots & \vdots \\ \vdots & & & \ddots & \vdots & \vdots \\ \varepsilon_p & \cdots & \cdots & \varepsilon_p & 1-(p-1)\varepsilon_p \end{pmatrix}$$

In this case $F$ is a map of $[-1, 1]^p$ into itself.

## 4.4. Two and Three Very Very Weakly Coupled Symmetric Tent Maps with Double Precision Number

In this section, only symmetric tent maps with $a = 2$ are considered. The computations are done with the same Dell computer and the same Borland C compiler as previously. Double precision numbers are used. First we verify for several values of $\varepsilon_1$ that for the 2-coupled maps case, each component $x_n$ and $y_n$ of the iterates of (4.6) are equally distributed over the interval $[-1, 1]$. That is, the invariant measure of each component converges to the Lebesgue measure.

### 4.4.1. Approximated Distribution Function

In order to compute numerically an approximation of the invariant measure also called the probability distribution function $P_N(x)$ (see section 3.2.1) linked to the one dimensional map $f$, we build a regular partition of $M$ small intervals (boxes) of the considered interval $[-1, 1]$:

$$r_i = [s_i, s_{i+1}[ \, , \, i = 0, M - 2 \quad (4.12)$$

$$r_{M-1} = [s_{M-1}, 1] \quad (4.13)$$

$$s_i = -1 + \frac{2i}{M} \quad i = 0, M \quad (4.14)$$

the length of which is:

$$s_{i+1} - s_i = \frac{2}{M} \quad (4.15)$$

We collect all iterates $f^{(n)}(x)$ belonging to these boxes (after a transient regime of $q$ iterations decided a priori, i.e. the first $q$ iterates are neglected). Once the computation of $N + q$ iterates is completed, the relative number of iterates with respect to $N/M$ in each box $r_i$ represents the value $P_N(s_i)$. The approximated probability distribution function $P_N(x)$ defined in this article is then a step function, with $M$ steps. As $M$ can vary in the next sections, we define :

$$P_{M,N}(s_i) = \frac{1}{2}\frac{M}{N}(\# r_i) \tag{4.16}$$

where $\# r_i$ is the number of iterates belonging to the interval $r_i$ and the constant $\frac{1}{2}$ allows the normalisation of $P_{M,N}(x)$ on the interval [-1, 1].

$$P_{M,N}(x) = P_{M,N}(s_i) \quad \forall x \in r_i \tag{4.17}$$

In the case of coupled maps, we are interested by the distribution of the component $x^1$, ..., $x^p$ of $X$ rather than the distribution of variable $X$ itself in $[-1, 1]^p$. Then we will consider the approximated probability distribution function $P_N(x^j)$ associated to one or several components of $F(X)$ defined by (4.7) which are one-dimensional maps.

*4.4.2. Distribution of Iterates of 2- and 3-coupled Symmetric Tent Maps*

As intuitively expected, Fig. 14 shows the convergence of the density of iterates of the components of 2-coupled symmetric tent maps to the Lebesgue measure when $\varepsilon_1$ converges towards 0. We call also discretization number ($N_{disc}$) the number $M$ of boxes of the interval. In the numerical experiments reported here, this number is equal to 100,000. We use the notation $N_{iter}$ (number of iterates) for $N$. In these experiments $N_{iter}$ is fixed to one hundred millions.

Moreover, for a fixed value of $N_{disc}$ when the number $N_{iter}$ increases, the discrepancy between $P_{N_{disc}, N_{iter}}(x)$ and the Lebesgue measure is expected to converge towards 0, except if there exist periodic orbits of finite length lower than $N_{iter}$ which capture the iterates. In this case whatsoever the value of $N_{iter}$ is, the approximated distribution function converges to the distribution function of the periodic orbit if it is unique or to the average of the distribution functions of the periodic orbits observed if not.

The discrepancies $E_1$ (i.e. in norm $L_1$) and $E_2$ (in norm $L_2$) between $P_{N_{disc}, N_{iter}}(x)$ and the Lebesgue measure are defined by:

$$E_{1,x}(N_{disc}, N_{iter}) = \left\| P_{N_{disc}, N_{iter}}(x) - 0.5 \right\|_{L_1} \tag{4.18}$$

$$E_{2,x}(N_{disc}, N_{iter}) = \left\| P_{N_{disc}, N_{iter}}(x) - 0.5 \right\|_{L_2} \tag{4.19}$$

Figs. 15 and 16 show the errors $E_{1,x}(N_{disc}, N_{iter})$ and $E^2_{2,x}(N_{disc}, N_{iter})$ versus the number of iterates of the approximated distribution functions with respect to the first variable $x$ for 2- and 3-coupled symmetric tent map. $N_{disc}$ is fixed to $10^5$, $\varepsilon_1$ to $10^{-14}$, $N_{iter}$ varies from $10^5$ to $5.10^{11}$ for the 2-coupled case and to $10^{12}$ for the 3-coupled one. Numerical results are given in Table 10 for 2-coupled maps and in Table 11 for 3-coupled maps[8].

---

[8] **Remark.** In order to made easier the comparison of the results, we display the square of the discrepancy $E_2$ instead of $E_2$ itself, the discrepancy being divided by 10 each time the number of iterations is multiplied by 10.

Same results are obtained for the variables *y* or *z*.

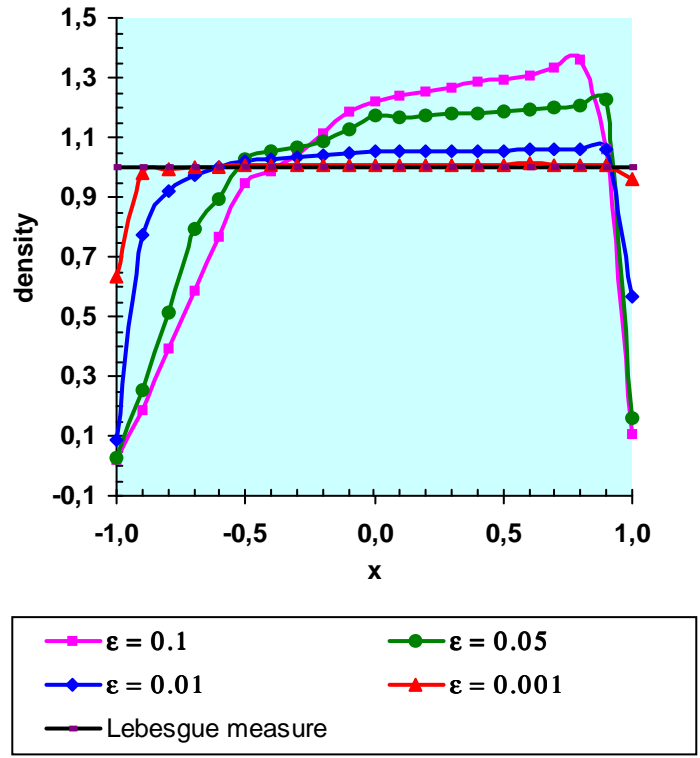

**Fig. 14.** Density of iterates of 2-coupled symmetric tent maps, double precision, $N_{discretization} = 10^5$, $\varepsilon_i = i\,\varepsilon_1$, $\varepsilon_1 = 10^{-1}$ to $10^{-3}$, $N_{iter} = 10^8$. Initial values $x_0 = 0.330$, $y_0 = 0.3387564$.

**Table 10.** Errors $L_1$ and $(L_2)^2$ for 2-coupled symmetric tent maps: Computations using double precision numbers (14 - 15 digits), $\varepsilon_i = i\,\varepsilon_1$, initial values $x^1_0 = 0.330$, $x^2_0 = 0.3387564$.

| $\varepsilon_1$ | $N_{disc}$ | $N_{iter}$ | $E_{1,x}(N_{disc}, N_{iter})$ | $E^2_{2,x}(N_{disc}, N_{iter})$ |
|---|---|---|---|---|
| $10^{-14}$ | $10^5$ | $10^5$ | $7.393 \times 10^{-1}$ | $1.015$ |
| $10^{-14}$ | $10^5$ | $10^6$ | $2.506 \times 10^{-1}$ | $1.006 \times 10^{-1}$ |
| $10^{-14}$ | $10^5$ | $10^7$ | $7.999 \times 10^{-2}$ | $1.004 \times 10^{-2}$ |
| $10^{-14}$ | $10^5$ | $10^8$ | $2.528 \times 10^{-2}$ | $1.004 \times 10^{-3}$ |
| $10^{-14}$ | $10^5$ | $10^9$ | $8.237 \times 10^{-3}$ | $1.066 \times 10^{-4}$ |
| $10^{-14}$ | $10^5$ | $10^{10}$ | $3.338 \times 10^{-3}$ | $1.799 \times 10^{-5}$ |
| $10^{-14}$ | $10^5$ | $10^{11}$ | $2.415 \times 10^{-3}$ | $9.189 \times 10^{-6}$ |
| $10^{-14}$ | $10^5$ | $5.10^{11}$ | $2.314 \times 10^{-3}$ | $8.421 \times 10^{-6}$ |

In both cases no periodic orbit is found even if the computation is completed to one trillion of iterations. In comparison with the Lozi map where one **G**-periodic orbit is found computing less than seven hundreds billions of iterates, it seems that the length of the periodic orbit of 3-coupled maps should be greater.

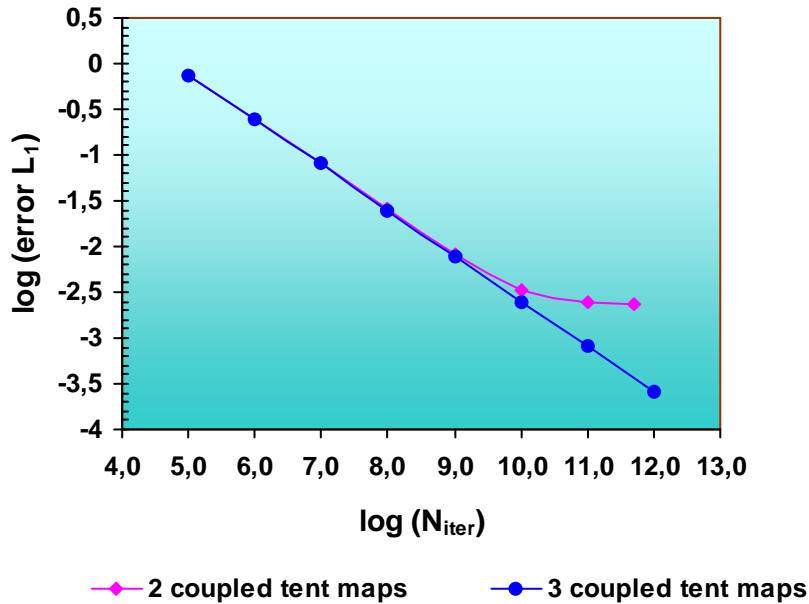

**Fig. 15.** Error $L_1$ for 2 and 3-coupled symmetric tent map, double precision, $N_{discretization} = 10^5$, $\varepsilon_i = i\, \varepsilon_1$, $\varepsilon_1 = 10^{-14}$, $N_{iter} = 10^5$ to $10^{12}$.
Initial values $x^1_0 = 0.330$, $x^2_0 = 0.3387564$, $x^3_0 = 0.331353429$.

3-coupled maps : correlation coefficient of $\dfrac{\log_{10}(errorL_1)}{\log_{10}(N_{iter})} = -0.999946$.

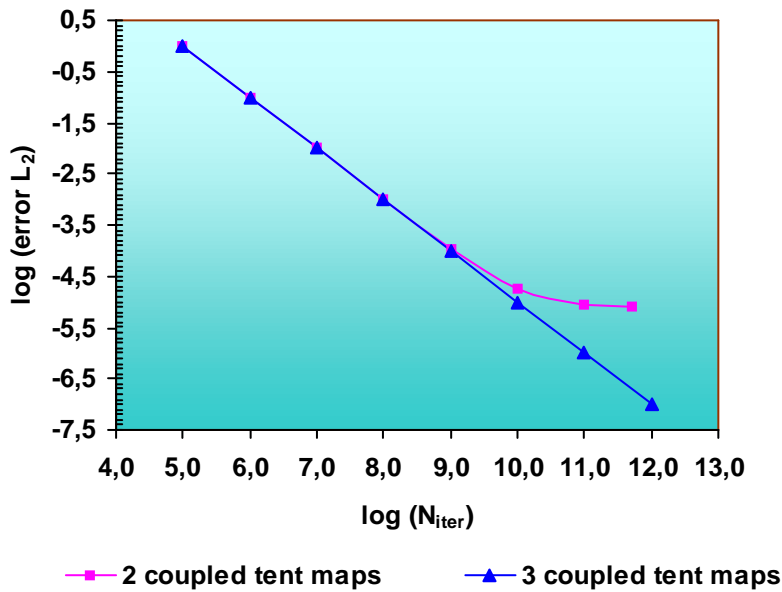

**Fig. 16.** (Error $L_2)^2$ for 2 and 3-coupled symmetric tent map, double precision, $N_{discretization} = 10^5$, $\varepsilon_i = i\, \varepsilon_1$, $\varepsilon_1 = 10^{-14}$, $N_{iter} = 10^5$ to $10^{12}$.
Initial values $x^1_0 = 0.330$, $x^2_0 = 0.3387564$, $x^3_0 = 0.331353429$.

3-coupled maps : correlation coefficient of $\dfrac{\log_{10}(errorL_2)^2}{\log_{10}(N_{iter})} = -0.999997$.

**Table 11.** Errors $L_1$ and $(L_2)^2$ for 3-coupled symmetric tent maps: Computations using double precision numbers (14 - 15 digits), $\varepsilon_i = i\,\varepsilon_1$, initial values $x^1_0 = 0.330$, $x^2_0 = 0.3387564$, $x^3_0 = 0.33153429$.

| $\varepsilon_1$ | $N_{discretization}$ | $N_{iter}$ | $E_{1,x}(N_{disc}, N_{iter})$ | $E^2_{2,x}(N_{disc}, N_{iter})$ |
|---|---|---|---|---|
| $10^{-14}$ | $10^5$ | $10^5$ | $7.383 \times 10^{-1}$ | $1.012$ |
| $10^{-14}$ | $10^5$ | $10^6$ | $2.497 \times 10^{-1}$ | $9.985 \times 10^{-2}$ |
| $10^{-14}$ | $10^5$ | $10^7$ | $8.000 \times 10^{-2}$ | $1.002 \times 10^{-2}$ |
| $10^{-14}$ | $10^5$ | $10^8$ | $2.518 \times 10^{-2}$ | $9.945 \times 10^{-4}$ |
| $10^{-14}$ | $10^5$ | $10^9$ | $7.910 \times 10^{-3}$ | $9.851 \times 10^{-5}$ |
| $10^{-14}$ | $10^5$ | $10^{10}$ | $2.495 \times 10^{-3}$ | $9.769 \times 10^{-6}$ |
| $10^{-14}$ | $10^5$ | $10^{11}$ | $8.011 \times 10^{-4}$ | $1.007 \times 10^{-6}$ |
| $10^{-14}$ | $10^5$ | $10^{12}$ | $2.541 \times 10^{-4}$ | $1.014 \times 10^{-7}$ |

The curves plotted on Figs. 15 and 16 show two characteristics:

1. For the 2-coupled map $E_{1,x}(N_{disc}, N_{iter})$ is decreasing to $2.3 \times 10^{-3}$ and $E^2_{2,x}(N_{disc}, N_{iter})$ is decreasing to $8.4 \times 10^{-6}$, however the convergence of both errors seems lower bounded by a minimal error.

2. Instead for the 3-coupled maps there is no such lower bound since the curves of the logarithms of the errors $L_1$ and $(L_2)^2$ versus the number of iterates look as straight lines in the range of the computations. It is even possible to compute the correlation coefficient which strengthen this belief:

$$\text{Correlation coefficient of } \frac{\log_{10}(E_{1,x}(N_{disc}, N_{iter}))}{\log_{10}(N_{iter})} = -0.999946$$

$$\text{Correlation coefficient of } \frac{\log_{10}(E^2_{2,x}(N_{disc}, N_{iter}))}{\log_{10}(N_{iter})} = -0.999997$$

For $N_{disc} = 100,000$ the numerically founded equations of these straight lines are:

$$E_{1,x}(N_{disc}, N_{iter}) = 10^{(-0.496821 \cdot \log_{10}(N_{iter}) + 2.369699)} \quad (4.20)$$

$$E^2_{2,x}(N_{disc}, N_{iter}) = 10^{(-1.000151 \cdot \log_{10}(N_{iter}) + 5.000686)} \quad (4.21)$$

In conclusion, the 3-coupled symmetric tent maps model we propose with very very small value of $\varepsilon_1$, seems a sterling model of generator of chaotic numbers with a uniform distribution of these numbers over the interval [-1, 1].

In order to go into detail for the study of its properties, the computations of more iterates is necessary. However, as it seems that no irregularity is found with up to one trillion of iterates computed using double precision numbers and in order to avoid computations of more than one or two weeks for every test on a desktop computer, we turn now our investigations to floating numbers instead of double precision ones.

### 4.5. Multiple Weakly Coupled Symmetric Tent Maps Computed with Floating Numbers

In this section, we will study the numerical properties of model (4.8 - 4.10) varying successively the number of coupled maps, the coupling constant $\varepsilon_1$, the number of iterations $N_{iter}$, and the discretization of the interval $N_{disc}$.

*4.5.1. Variation of the Number of Coupled Maps*

We consider now numerical experiments always done with the same Dell computer and the same Borland C compiler. Instead of double precision numbers, simple precision numbers or floating numbers (i.e. with 7-8 digits) are used. We compare the discrepancy between $P_{N_{disc}, N_{iter}}(x)$ and the Lebesgue measure for $p$-coupled maps, $p$ being equal to 2 up to 7. The coupling constant $\varepsilon_1$ is fixed to $10^{-7}$, $N_{disc}$ is $10^5$ and $N_{iter}$ varies from $10^5$ up to $10^{11}$.

The results are given on Table 12.

**Table 12.** $E_{1,x}(N_{disc}, N_{iter})$ **for n-coupled symmetric tent maps: Computations using simple precision numbers (7-8 digits) with respect to $N_{iter}$.**
$N_{disc} = 10^5$, $N_{iter} = 10^5$ to $10^{11}$, $\varepsilon_1 = 10^{-7}$, $\varepsilon_i = i\,\varepsilon_1$.
Initial values $x^1_0 = 0.330$, $x^2_0 = 0.3387564$, $x^3_0 = 0.3313534$, $x^4_0 = 0.3332135$, $x^5_0 = 0.3387325$, $x^6_0 = 0.3438542$, $x^7_0 = 0.3654218$.

| $N_{iter}$ | 2-maps | 3-maps | 4-maps | 5-maps | 7-maps |
|---|---|---|---|---|---|
| $10^5$ | | $7.302 \times 10^{-1}$ | $7.388 \times 10^{-1}$ | $7.356 \times 10^{-1}$ | $7.345 \times 10^{-1}$ |
| $10^6$ | $3.219 \times 10^{-1}$ | $2.493 \times 10^{-1}$ | $2.517 \times 10^{-1}$ | $2.502 \times 10^{-1}$ | $2.489 \times 10^{-1}$ |
| $10^7$ | $2.929 \times 10^{-1}$ | $8.004 \times 10^{-2}$ | $7.949 \times 10^{-2}$ | $7.967 \times 10^{-2}$ | $7.987 \times 10^{-2}$ |
| $10^8$ | $2.964 \times 10^{-1}$ | $2.581 \times 10^{-2}$ | $2.531 \times 10^{-2}$ | $2.522 \times 10^{-2}$ | $2.530 \times 10^{-2}$ |
| $10^9$ | | $9.362 \times 10^{-3}$ | $8.353 \times 10^{-3}$ | $8.102 \times 10^{-3}$ | $8.039 \times 10^{-3}$ |
| $10^{10}$ | | $5.648 \times 10^{-3}$ | $3.501 \times 10^{-3}$ | $2.697 \times 10^{-3}$ | $2.682 \times 10^{-3}$ |
| $10^{11}$ | | $5.159 \times 10^{-3}$ | $2.510 \times 10^{-3}$ | $1.233 \times 10^{-3}$ | $1.198 \times 10^{-3}$ |
| periodic solution | period 1,320,572 *Megaperiodic solution* | period 76,355,473,953 *Gigaperiodic solution* | Not found | Not found | Not found |

For 2-coupled maps one **M**-periodic orbit is found and one **G**-periodic orbit appears for 3-coupled maps. Instead computing up to one hundred billions of iterates no other periodic orbits are found for 4, 5 and 7-coupled maps. The discrepancies between the approximated invariant functions and the Lebesgue measure (logarithms of $E_{1,x}(N_{disc}, N_{iter})$) versus $N_{iter}$ for these values of $p$ are plotted on the Fig. 17. They are compared on the same figure to the results obtained in the previous subsection for 2 and 3-coupled maps using double precision numbers.

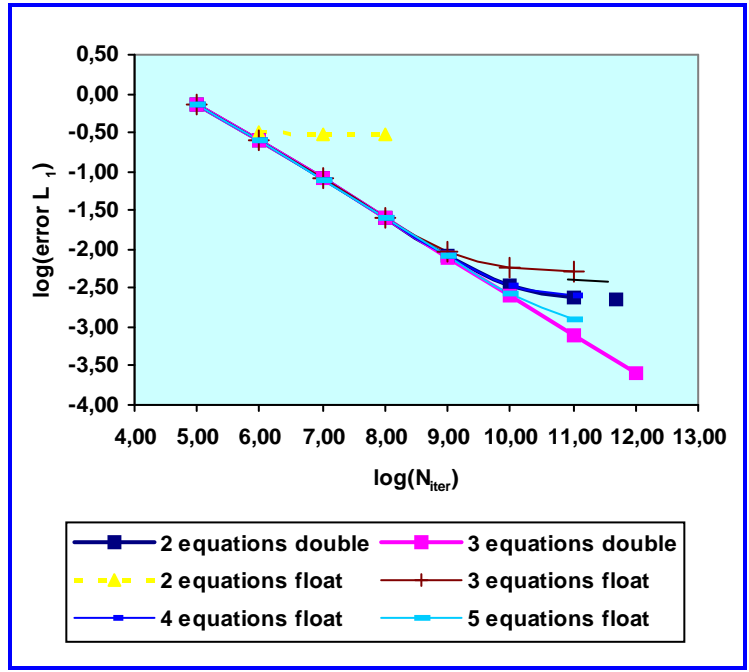

**Figure 17.** $E_{1,x}(N_{disc}, N_{iter})$ **for n-coupled symmetric tent maps using simple precision numbers (7-8 digits) and 2- and 3-coupled maps using double precision numbers, $N_{disc} = 10^5$, $\varepsilon = 10^{-7}$, $N_{iter} = 10^5$ to $10^{11}$. Initial values $x^1_0 = 0.330$, $x^2_0 = 0.3387564$, $x^3_0 = 0.3313534$, $x^4_0 = 0.3332135$, $x^5_0 = 0.3387325$, $x^6_0 = 0.3438542$, $x^7_0 = 0.3654218$.**

If we except the model of 2-coupled maps computed with floating numbers, all the models give good results for the production of chaotic numbers.

In the range of $10^5$ to $10^{11}$ iterates, the 4-coupled maps model with floating points gives results alike the 2-coupled maps models with double precision numbers.

It is not possible to distinguish between the results of the 5-coupled maps and the 7-coupled maps in this range of $N_{iter}$.

However, the 3-coupled maps model with double precision numbers is better than all the models with up to 7-coupled maps computed with floating numbers when $N_{iter} \geq 10^{11}$.

### 4.5.2. Variation of $\varepsilon_1$

We restrict ourselves to the 7-coupled maps model. The impact of variation of the coupling constant $\varepsilon_1$ is then considered. $N_{disc}$ is still fixed to $10^5$, and $N_{iter}$ to $10^{10}$.

The coupling constant $\varepsilon_1$ decreases from $10^{-2}$ down to $10^{-8}$. The results are given on Table 13 and plotted on Fig. 18.

The best values in the case of floating numbers are obtained when $\varepsilon_1$ belongs to the interval $[5 \cdot 10^{-8}, 10^{-7}]$.

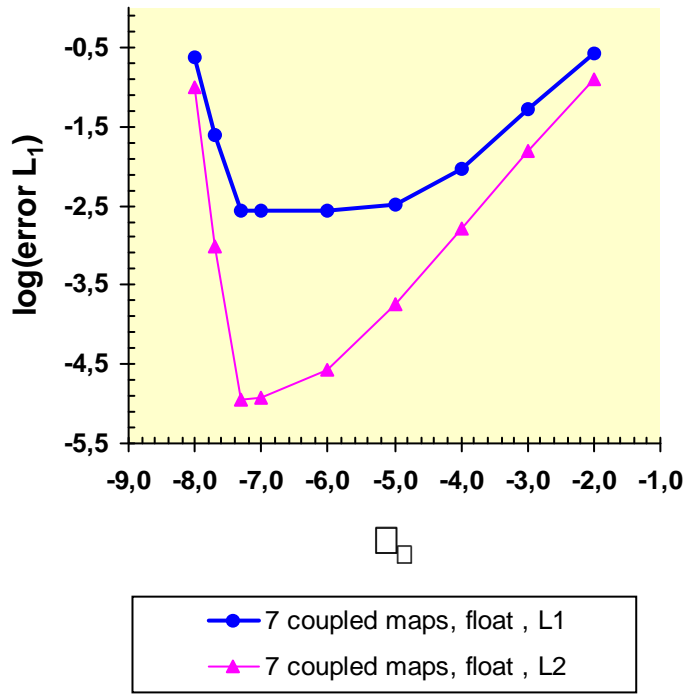

**Fig. 18.** Errors $L_1$ and $L_2$ for 7-coupled symmetric tent maps using simple precision numbers (7-8 digits) with respect to ε. $N_{\text{discretization}} = 10^5$, $N_{\text{iter}} = 10^{10}$, $ε_i = i\, ε_1$, $ε_1 = 10^{-8}$ to $10^{-2}$. Initial values $x^1_0 = 0.330$, $x^2_0 = 0.3387564$, $x^3_0 = 0.3313534$, $x^4_0 = 0.3332135$, $x^5_0 = 0.3387325$, $x^6_0 = 0.3438542$, $x^7_0 = 0.3654218$.

The precision of the floating numbers being between 7 and 8 digits, the best value of the coupling constant is near the value of their two last significant digits. In the case of double precision numbers similar results can be found. In practice we have chosen the value $10^{-14}$.

### 4.5.3. Variation of the Number of Iterates and the Discretization Number

We consider again the 7-coupled maps model. The impact of variation of the number of iterations $N_{\text{iter}}$ is then studied together with the variation of the number of boxes of the discretization $N_{\text{disc}}$.

We compare the discrepancy between $P_{N_{disc}, N_{iter}}(x)$ and the Lebesgue measure for this model. The coupling constant $ε_1$ is fixed to $10^{-7}$, $N_{\text{disc}}$ varies from $10^2$ to $10^5$, and $N_{\text{iter}}$ varies from $10^6$ up to $10^{11}$. The results are plotted on Figs. 19 and 20.

Behaviour of the error is very regular (except for the value $N_{\text{disc}} = 10^6$ and $N_{\text{iter}} = 10^{11}$ and in a less obvious manner for $N_{\text{disc}} = 10^6$ and $N_{\text{iter}} = 10^{10}$).

For example for $N_{\text{disc}} = 10^4$ the results are plotted on a straight line (brown straight line of Fig. 19 with triangular dots) with respect to the variable $N_{\text{iter}}$ the numerically founded equation of which is

$$E_{1,x}(N_{disc}, N_{iter}) = 10^{(-0.495261 \cdot \log_{10}(N_{iter}) + 1.869088)} \qquad (4.21)$$

**Table 13.** Errors $L_1$ and $(L_2)^2$ for 7-coupled symmetric tent maps: Computations using simple precision numbers (7-8 digits) with respect to $\varepsilon_1$.
$N_{disc} = 10^5$, $N_{iter} = 10^{10}$, $\varepsilon_1 = 10^{-8}$ to $10^{-2}$, $\varepsilon_i = i\,\varepsilon_1$. Initial values $x^1_0 = 0.330$, $x^2_0 = 0.3387564$, $x^3_0 = 0.3313534$, $x^4_0 = 0.3332135$, $x^5_0 = 0.3387325$, $x^6_0 = 0.3438542$, $x^7_0 = 0.3654218$.

| $\varepsilon_1$ | $E_{1,x}(N_{disc}, N_{iter})$ | $E_{2,x}(N_{disc}, N_{iter})$ |
|---|---|---|
| $10^{-2}$ | $2.616 \cdot 10^{-1}$ | $1.137 \cdot 10^{-1}$ |
| $10^{-3}$ | $5.343 \cdot 10^{-2}$ | $1.549 \cdot 10^{-2}$ |
| $10^{-4}$ | $8.980 \cdot 10^{-3}$ | $1.677 \cdot 10^{-3}$ |
| $10^{-5}$ | $3.209 \cdot 10^{-3}$ | $1.805 \cdot 10^{-4}$ |
| $10^{-6}$ | $2.717 \cdot 10^{-3}$ | $2.641 \cdot 10^{-5}$ |
| $10^{-7}$ | $2.682 \cdot 10^{-3}$ | $1.207 \cdot 10^{-5}$ |
| $5 \times 10^{-8}$ | $2.670 \cdot 10^{-3}$ | $1.146 \cdot 10^{-5}$ |
| $10^{-8}$ | $2.396 \cdot 10^{-1}$ | $9.816 \cdot 10^{-2}$ |

and for $N_{iter} = 10^{10}$ the results are plotted on a straight line (brown straight line of Fig. 20 with squared dots) with respect to $N_{disc}$ of which the numerically founded equation is

$$E_{1,x}(N_{disc}, N_{iter}) = 10^{(0.518417 \cdot \log_{10}(N_{disc}) - 5.140626)} \tag{4.22}$$

It is even possible to combine all the results in the following general equation of one plane versus both variables $N_{disc}$ and $N_{iter}$ as shown in Fig. 21 and to find numerically the fundamental equation

$$E_{1,x}(N_{disc}, N_{iter}) = 10^{(0.491388 \cdot (\log_{10}(N_{disc}) - (N_{iter})) - 0.104101)} \tag{4.23}$$

The correlation coefficient of the line (4.23) being $\dfrac{\log_{10}(E_{1,x}(N_{disc}, N_{iter}))}{\log_{10}(N_{disc} - N_{iter})} = 0.999081$

## 5. Very Weakly Coupled Symmetric Logistic Maps

We now consider again the logistic map (3.8) defined on the interval [-1,1]. We study the coupling of $p$ logistic maps following (4.10) and (4.11). Note that it is easy to produce equally distributed numbers on the interval, using coupled logistic map because the invariant measure (3.10) can be inverted explicitly.

### 5.1. Three and *p*-Coupled Symmetric Maps Computed with Floating Numbers

The numerical results obtained studying the logistic map instead of the symmetric tent map are similar but less powerful. For example for 2-coupled logistic maps, using simple precision numbers with a coupling constant $\varepsilon_1$ fixed to $10^{-7}$:

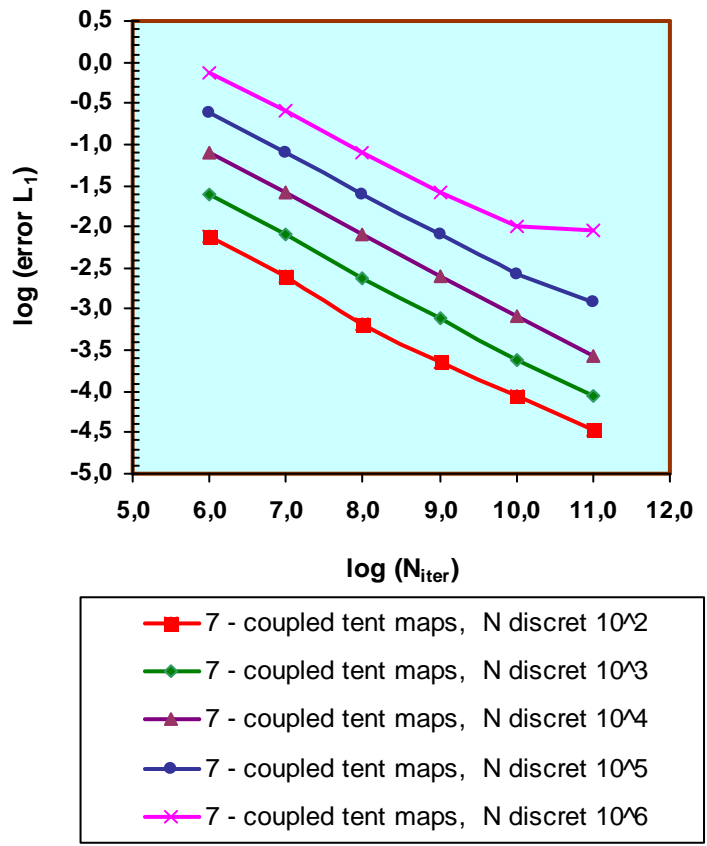

**Fig. 19.** Error $L_1$ for 7-coupled symmetric tent maps using simple precision numbers (7-8 digits) with respect to $N_{iter}$. $N_{disc} = 10^2$ to $10^6$, $\varepsilon_i = i\,\varepsilon_1$, $\varepsilon_1 = 10^{-7}$, $N_{iter} = 10^6$ to $10^{11}$. Initial values $x^1_0 = 0.330$, $x^2_0 = 0.3387564$, $x^3_0 = 0.3313534$, $x^4_0 = 0.3332135$, $x^5_0 = 0.3387325$, $x^6_0 = 0.3438542$, $x^7_0 = 0.3654218$.

if $\varepsilon_2 = 2\varepsilon_1$, a 477-periodic orbit is found (for the symmetric tent map a 1,320,572-periodic orbit is found when $\varepsilon_2 = 2\,\varepsilon_1$), if $\varepsilon_2 = 3\varepsilon_1$, a 16,917-periodic orbit is found, if $\varepsilon_2 = 0.5\varepsilon_1$, a 1,006,450-periodic orbit (**M**-periodic) is found.

For 3-coupled logistic maps, using simple precision numbers with a coupling constant $\varepsilon_1$ fixed to $10^{-7}$:

if $\varepsilon_i = i\varepsilon_1$, a 2,399,718,074-periodic orbit is found, (for the symmetric tent map a 76,355,473,953-periodic orbit is found with same conditions), if $\varepsilon_2 = 3\varepsilon_1$, a 16,917-periodic orbit is found, if $\varepsilon_2 = 0.5\varepsilon_1$, a 1,006,450-periodic orbit is found.

We will study more precisely this kind of results (dependence upon the variation of the coupling constants ratios) in a forthcoming paper.

In order to avoid numerical artifacts dues to the pike-shaped invariant measure (see Fig. 8) of the

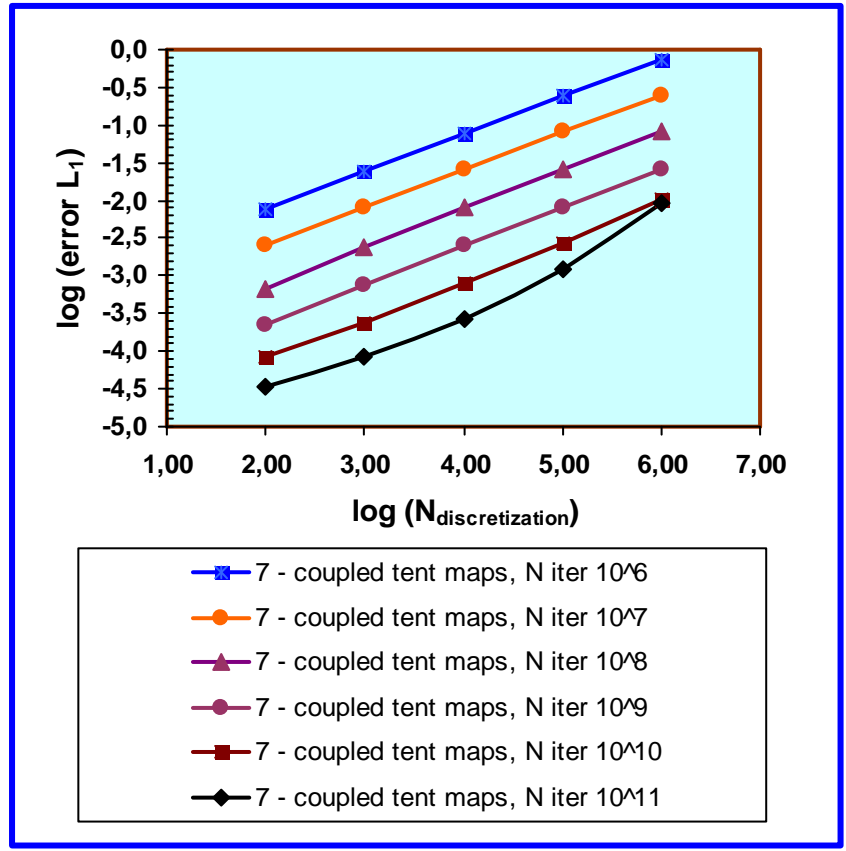

**Fig. 20.** Error $L_1$ for 7-coupled symmetric tent maps using simple precision numbers (7-8 digits) with respect to $N_{discretization}$. $N_{disc} = 10^2$ to $10^6$, $\varepsilon_i = i\,\varepsilon_1$, $\varepsilon_1 = 10^{-7}$, $N_{iter} = 10^6$ to $10^{11}$. Initial values $x^1_0 = 0.330$, $x^2_0 = 0.3387564$, $x^3_0 = 0.3313534$, $x^4_0 = 0.3332135$, $x^5_0 = 0.3387325$, $x^6_0 = 0.3438542$, $x^7_0 = 0.3654218$.

logistic map near the ends points of the interval [-1,1], we define

$$E^{**}_{1,x}(N_{disc}, N_{iter}) = \left|\ \|P_{N_{disc},N_{iter}}(x) - P(x)\|_{L_1}\bigg|_{[-0.98,-0.98]}\right.\quad (5.1)$$

with $P(x)$ given by (3.10), that is, we consider only the discrepancy between $P_{N_{disc},N_{iter}}(x)$ and $P(x)$ in the interval [-0.98, 0.98].

On Fig. 22 the errors $E^{**}_{1,x}(N_{disc}, N_{iter})$ versus the number of iterates with respect to the first variable $x$ for 3, 4, up to 7-coupled logistic map are plotted. $N_{disc}$ is fixed to $10^5$, $\varepsilon_1$ to $10^{-7}$, $N_{iter}$ varies from $2.10^9$ to $2.10^{10}$.

The same kind of results already obtained for the tent map in Fig. 17 are observed with the logistic map. The use of $n + 1$ coupled logistic maps is better than the use of $n$-coupled ones, however it seems that the process of improvement is not meaningful using more than 7-coupled maps.

For example, using quadratic interpolation of the black curve we obtain the value $E^{**}_{1,x}(N_{disc}, 10^{10}) = 0.0027552$ for the 7-coupled logistic maps and $E_{1,x}(N_{disc}, 10^{10}) = 0.0026816$ for the tent map (Table 12). These results are comparable, however in the case of the tent maps the computation is faster because no inversion of the invariant measure is needed. Moreover one can use the whole interval [-1, 1] instead of the reduced interval [-0.98, 0.98].

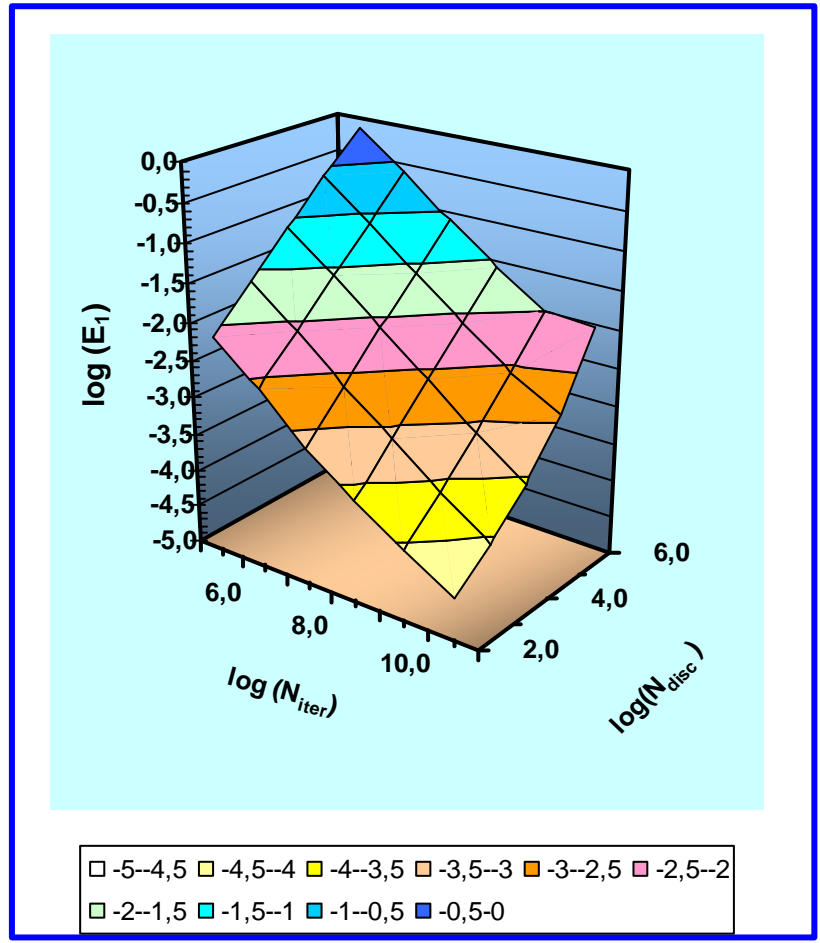

**Fig. 21.** Error $L_1$ for 7-coupled symmetric tent maps using simple precision numbers (7-8 digits) with respect to both $N_{iteration}$ and $N_{discretization}$. $N_{disc} = 10^2$ to $10^6$, $\varepsilon_i = i\,\varepsilon_1$, $\varepsilon_1 = 10^{-7}$, $N_{iter} = 10^6$ to $10^{11}$. Initial values $x^1_0 = 0.330$, $x^2_0 = 0.3387564$, $x^3_0 = 0.3313534$, $x^4_0 = 0.3332135$, $x^5_0 = 0.3387325$, $x^6_0 = 0.3438542$, $x^7_0 = 0.3654218$.

Same results are obtained for the variables *y* or *z*.

**5.2. More Chaotic Numbers**

As said before, in the case of coupled maps, we are interested by the distribution of the component $x^1, \ldots, x^p$ of $X$ rather than the distribution of variable $X$ itself in $[-1, 1]^p$. Then we will consider the approximated probability distribution function $P_N(x^j)$ associated to one or several components of $F(X)$ defined by (4.7) which are one-dimensional maps. However, we compute in

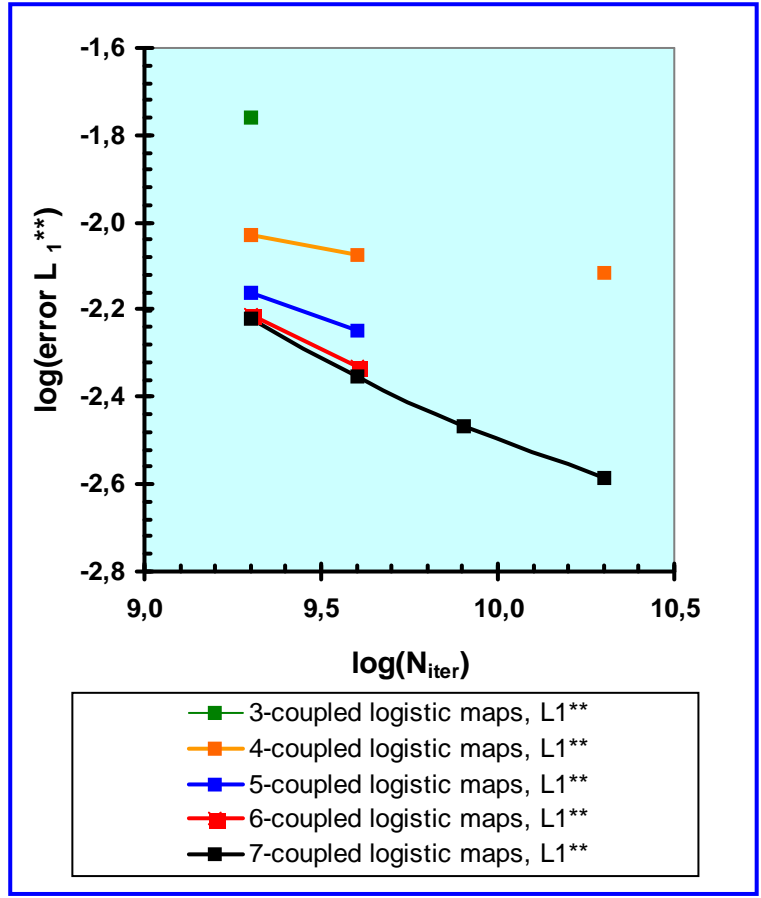

**Fig. 22.** $E_{1,x}^{**}(N_{disc}, N_{iter})$ **for 3, 4, 5, 6 & 7-coupled symmetric logistic maps using simple precision numbers (7-8 digits), $N_{disc} = 10^5$, $\varepsilon = 10^{-7}$, $N_{iter} = 2.10^9$ to $2.10^{10}$. Initial values $x^1_0 = 0.330$, $x^2_0 = 0.3387564$, $x^3_0 = 0.3313534$, $x^4_0 = 0.3332135$, $x^5_0 = 0.3387325$, $x^6_0 = 0.3438542$, $x^7_0 = 0.3654218$.**

fact simultaneously $p$ sequences of chaotic numbers. It could be interesting to use all these sequences together in order to produce chaotic number $p$ times faster. One can mix these sequences in the simple manner $\{\underline{x}_n\} = \ldots, x^1_n, x^2_n, \ldots, x^p_n, x^1_{n+1}, x^2_{n+1}, \ldots, x^p_{n+1}, x^1_{n+2}, \ldots$. In Fig. 23, the errors $E_{1,x}(N_{disc}, N_{iter})$ and $E_{1,x}^{**}(N_{disc}, N_{iter})$ for 6 and 7-coupled logistic maps are plotted. There is a slight improvement of the convergence of the approximated probability distribution function $P_N(\underline{x})$ associated to the mixed sequence $\{\underline{x}_n\}$ with respect to the approximated probability distribution function $P_N(x^j)$ associated to one component (we have plotted the errors $E_{1,\underline{x}}(N_{disc}, N_{iter})$ and $E_{1,\underline{x}}^{**}(N_{disc}, N_{iter})$ taken into account the actual number of chaotic numbers which is equal to $p \times N_{iter}$)

This shows that the use of the mixed sequence of chaotic numbers $\{\underline{x}_n\}$ is very efficient.

The interpretation of the slight difference between the results obtained with 6 or 7 coupled maps is however unclear at this time and will need more computations.

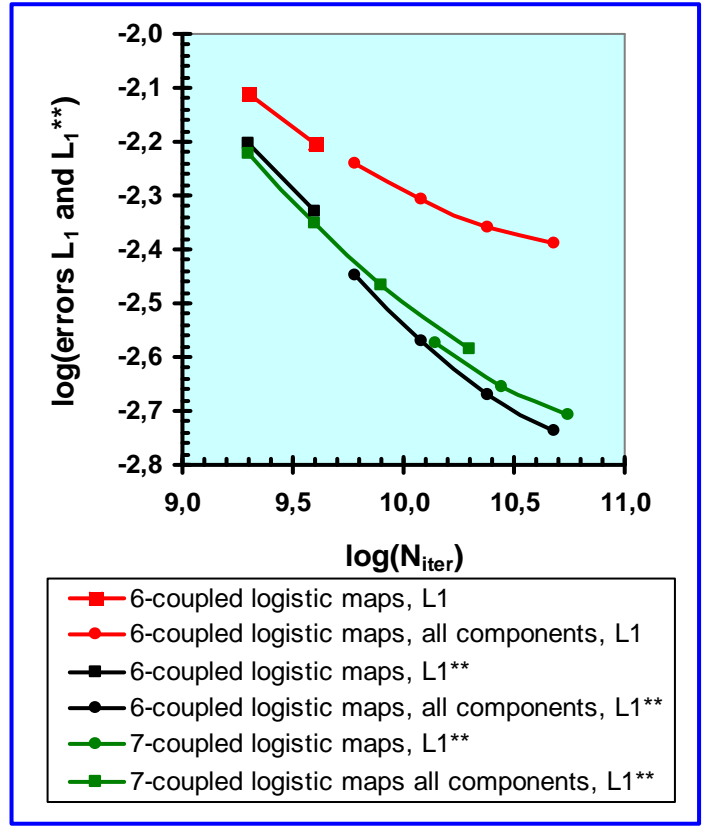

**Fig. 23.** $E_{1,x}^{**}(N_{disc}, N_{iter})$ for 6 & 7-coupled symmetric logistic maps using simple precision numbers (7-8 digits), $N_{disc} = 10^5$, $\varepsilon = 10^{-7}$, $N_{iter} = 2.10^9$ to $6.10^{10}$. Error for one or all the components. Initial values $x^1_0 = 0.330$, $x^2_0 = 0.3387564$, $x^3_0 = 0.3313534$, $x^4_0 = 0.3332135$, $x^5_0 = 0.3387325$, $x^6_0 = 0.3438542$, $x^7_0 = 0.3654218$.

## 6. Conclusion

In conclusion, we have considered the very very weakly coupling of symmetric tent and logistic maps using single or double precision numbers. In particular the 3-coupled symmetric tent maps model we propose with very very small value of $\varepsilon_1$, seems a sterling model of generator of chaotic numbers with a uniform distribution of these numbers over the interval [-1, 1]. It produces also very long periodic orbits: Gigaperiodic orbits when computed with simple precision numbers, and orbits of unknown length when computed with double precision numbers. These chaotic sequences are not random sequences, however they can be used when chaos has been found more useful that random [4].

Modern Mathematical Models, Methods and Algorithms for Real World Systems, A.H. Siddiqi, I.S. Duff and O. Christensen (Editors), Anamaya Publishers, New Delhi, India, 80-14, 2006. *Proc. Conf. Intern. on Industrial and Appl. Math., New Delhi, India 4-6 Dec. 2004. Invited conference.*## References